\documentclass[12pt]{article}
\usepackage{amsthm}
\usepackage{latexsym}
\usepackage[utf8]{inputenc}
\usepackage{amsmath}
\usepackage{amssymb}
\usepackage{amsfonts}
\usepackage[T1]{fontenc} 
\usepackage{amsmath}
\usepackage{amsfonts}
\usepackage{amsfonts}
\usepackage{amsthm}
\usepackage{graphicx}
\usepackage{latexsym}
\makeatletter
\renewcommand{\@seccntformat}[1]
{\csname the#1\endcsname.\enspace}
\makeatother
\input{amssym.def}
\input{amssym}

\begin{document}
\newtheorem{theo}{Th\'eor\`eme}
\newtheorem{coro}{Corollaire}
\newtheorem{lem}{Lemme}
\newtheorem{definition}{Definition}[section]
\newtheorem{theorem}{Theorem}
\newtheorem{lemma}{Lemma}
\newtheorem{prop}{proposition}
\newtheorem{remark}{Remark}
\newtheorem{re}{Remark}
\newtheorem{corollary}{Corollary}
\newtheorem{example}{Example}
\renewcommand{\theequation}{\thesection.\arabic{equation}}
\renewcommand{\thetheorem}{\thesection.\arabic{theorem}}
\renewcommand{\thelemma}{\thesection.\arabic{lemma}}
\renewcommand{\thecorollary}{\thesection.\arabic{corollary}}
\renewcommand{\theremark}{\thesection.\arabic{remark}}
\renewcommand{\theexample}{\thesection.\arabic{example}}
\begin{center}
\textbf{Bayesian estimation and prediction for certain type of mixtures} 
\end{center}

\begin{center}
{\sc Aziz LMoudden$^{a,b}$ \& \'Eric Marchand$^{a}$, \footnote{\today}} \\
{\it a  Universit\'e de
    Sherbrooke, D\'epartement de math\'ematiques, Sherbrooke Qc,
    CANADA, J1K 2R1 \\
    b  Département de mathématiques, Université du Québec à Montréal, Montréal, Qc, CANADA, H2X 3Y7 \quad e-mails:  eric.marchand@usherbrooke.ca;\,\, aziz.lmoudden@usherbrooke.ca } \\
\end{center}
\begin{center}
{\sc Summary} 
\end{center}

\small  
\noindent  For two vast families of mixture distributions and a given prior, we provide unified representations of posterior and predictive distributions.   Model applications presented include bivariate mixtures of  Gamma distributions labelled as Kibble-type, non-central Chi-square and F distributions, the distribution of $R^2$ in multiple regression, variance mixture of normal distributions, and mixtures of location-scale exponential distributions including the multivariate Lomax distribution.   An emphasis is also placed on analytical representations and the relationships with a host of existing distributions and several hypergeometric functions of one or two variables.



\normalsize

\section{Introduction}

\noindent Mixture models are ubiquitous in probability and statistics.  Such models, whether they are finite mixture models, mixtures of Poisson, exponential, gamma or normal distributions, etc.,  are quite useful and appealing for best representing data and heterogeneous environments.  As well, distributional properties of mixture models are often quite elegant and instructive.    However, analytical and numerical challenges are present and well documented, namely in terms of likelihood-based and Bayesian inference.  

\noindent  It is also the case that several familiar distributions are representable in terms of mixtures, and that such representation facilitates the derivation of various statistical properties and approaches to statistical inference.   Prominent examples include the noncentral chi-square, Beta, and Fisher distributions, which arise typically in relationship to quadratic forms in normal linear models.  Other important examples include the distribution of the square of a multiple correlation coefficient ($R^2$) in a standard multiple regression linear model with normally distributed errors, as well as the vast class of univariate or bivariate Gamma mixtures which include the Kibble distribution (see Example \ref{examplekibble}).

\noindent We consider mixture models for summary statistics $X,Y \in \mathbb{R}^d$ which are of one of the following two types:
\begin{eqnarray}
\label{modeltypeI}
 &\hbox{Type I}: \,\,X|K \sim f_K \hbox{ with }  K|\theta \sim g_{\theta},\, & Y|J \sim q_J \hbox{ with }  J|\theta \sim h_{\theta}\,; \\ 
 \label{modeltypeII}
&\hbox{Type II}: \,\,X|K,\theta \sim f_{\theta,K} \hbox{ with }  K \sim g , & Y|J, \theta \sim q_{\theta,J} \hbox{ with }  J \sim h \,;
\end{eqnarray}
\noindent The classification of these types will be adhered to and seems to be a natural way to present the various expressions and examples that make up this paper.   In the above models, the mixing variables $K,J$ will typically be either a discrete or continuous univariate distribution, but more generally $K,J\in \cal{K}$ with $\cal{K}$ a measurable space.   The parameter $\theta\in \mathbb{R}^p$ is unknown and a prior density will be assumed.  
Otherwise the densities $f_k, f_{\theta,k}, q_j, \hbox{ and } q_{J,\theta}$ will be assumed known (except for the value of $\theta$) and absolutely continuous with respect to a $\sigma-$finite measure $\nu$.  Similarly, the densities $g$ and $g_{\theta}$ for $K$, as well as $h$ and $h_{\theta}$ for $J$, will be assumed to be known and absolutely continuous with respect to a $\sigma-$finite measure $\mu$.    Examples will include both discrete and continuous mixing for $K$ and $J$, as well as discrete and continuous models for the conditional distributions of $X$ and $Y$.

\noindent  With analytical expressions (Theorem \ref{theoremtypeI} and Theorem \ref{theoremtypeII}) for Bayesian posterior distributions for $\theta$ based on $X$, as well as Bayesian predictive densities for $Y$ based on $X$, we are particularly interested in eliciting the general structures driving these Bayesian solutions and illustrating the wide breadth of applications.  In the illustrations, the solutions involve mixtures themselves and bring into play special functions such as Hypergeometric and Appell which are readily available in mathematical packages.   They are also ``exact' and do not rely on stochastic approximations (e.g., MCMC) which are omnipresent in Bayesian computing.    

\noindent Much focus is placed on identifying concise representations and informative connections between distributions.  Model applications that we present include bivariate mixtures of  Gamma distributions that we label as Kibble-type, non-central Chi-square and F distributions, the distribution of $R^2$ in multiple regression, variance mixture of normal distributions, and mixtures of location-scale exponential distributions including the multivariate Lomax distribution.     For bivariate gamma mixtures, which we focus on in Section 4, we also consider bivariate prior distributions with dependence structures, in particular as occurring under an order restriction on the parameters.
The posterior and predictive distribution decompositions in this work follow familiar paths, but the analytical representations provided are nevertheless unified, useful and insightful, and lead to simplifications in Bayesian posterior analyses.  This is also the case where we purposely exploit the mixture representation of a familiar distribution.

\noindent  Here is a first illustration of Type I and Type II mixtures, which will be studied more generally in Sections 3 and 4 respectively, and which borrows notations presented in Section 2.

\begin{example}
\label{examplekibble}
The Kibble bivariate gamma distribution (Wicksell, 1933; Kibble, 1941) for $X=(X_1, X_2)$,  admits the following mixture representation:
 \begin{equation}
 \label{kibblerepresentation}
 X_i|K =^{\!\!\!\!^{i.d}} G(\nu+K, \frac{1-\rho}{\lambda_i}), i=1,2,  \hbox{ with }  K =^{\!\!\!^d} \hbox{NB}(\nu, 1-\rho)\,.
 \end{equation}

Here $\nu, \lambda_1, \lambda_2>0$, and $\rho \in [0,1).$ Such a distribution originates for instance in describing the joint distribution of sample variances generated from a bivariate normal distribution with correlation coefficient $\rho$.  The cases $\rho=0$ reduce to $X_i= ^{\!\!\!\!^{i.d}} \hbox{G}(\nu, \frac{1}{\lambda_i}) , i=1,2$.  Now, observe that we have a Type I mixture if $\theta=\nu$ with known $\lambda_1, \lambda_2, \rho$, and a Type II mixture for $\theta=(\lambda_1, \lambda_2)$; or again $\theta=(\frac{1-\rho}{\lambda_1},\frac{1-\rho}{\lambda_2} )$; with known $\nu, \rho$.  A third type of mixture arises when $\rho$ is unknown, but will not be further addressed in this paper.  The mixing representation of the Kibble distribution has been exploited in a critical fashion for Bayesian analysis by Iliopoulos et al. (2005).
\end{example}

\section{Notations and definitions}

\noindent  Here are some notations and definitions used throughout concerning some special functions and various distributions that appear below, related either to the model, the mixing variables, the prior, the posterior, or the predictive distribution.  We will write $X \sim f$ to indicate that the random variable $X$ has density (or p.m.f.) $f$.  We will also write $f \sim D$ to signify that density $f$ is associated with distribution $D$.   We will the use the symbol $\,=^{\!\!\!^d}\,$ to represent equality in distribution and $\,=^{\!\!\!\!^{i.d}}\,$ for equality in distribution with independent components.

\noindent In what follows, we denote $\hbox{P}(\theta)$ as the Poisson distribution with mean $\theta$ and p.m.f. (or density) $\frac{e^{-\theta} \theta^n}{n!} \, \mathbb{I}_{\mathbb{N}}(n)$.  We write $N =^{\!\!\!^d} NB(r,\xi)$, with $r>0$ and $\xi \in (0,1)$, to denote a negative binomial distribution with p.m.f. (or density) $\frac{(r)_n}{n!} \, \xi^r (1-\xi)^n  \, \mathbb{I}_{\mathbb{N}}(n)$. 

\noindent   Throughout the paper, we define for positive real numbers $a_1, \ldots, a_p, b_1, \ldots, b_q$, and $z \in \mathbb{R}$, the generalized hypergeometric function as
\begin{equation}
\nonumber
  _pF_q(a_1,...,a_p; b_1,...,b_q; z)=\sum_{k=0}^{\infty}\frac{\prod_{i=1}^p(a_i)_k}{\prod_{j=1}^q(b_j)_k}\frac{z^k}{k!}\,,
\end{equation}
with  $(\alpha)_m = \frac{\Gamma(\alpha+m)}{\Gamma(\alpha)}$ the Pochhammer function defined here for $\alpha>0, m \in \mathbb{N}$.  We write $N =^{\!\!\!^d} \hbox{Hyp}(a_1, \ldots a_p; b_1, \ldots, b_q; \lambda)$, for $a_i>0, b_j>0$, to denote a generalized hypergeometric distribution (e.g., Johnson et al., 1995) with p.m.f.
\begin{equation}
\nonumber  \hbox{P}(N=n) \, = \, \frac{1}{_pF_q(a_1,...,a_p; b_1,...,b_q; \lambda)} \,
\frac{\prod_{i=1}^p(a_i)_n}{\prod_{j=1}^q(b_j)_k}\frac{\lambda^n}{n!}\,\, \mathbb{I}_{\mathbb{N}}(n)\,.
\end{equation}
For these distributions,  we will have positive $a_i$'s and $b_j$'s, and $\lambda$ will be non negative and in the radius of convergence of the $_pF_q$ function.   With such a notation for instance, we could write $X =^{\!\!\!^d} \hbox{Hyp}(-;-;\theta)$ for $X =^{\!\!\!^d} \hbox{P}(\theta)$, and $X =^{\!\!\!^d} \hbox{Hyp}(r;-;1-\xi)$ for $X =^{\!\!\!^d} \hbox{NB}(r,\xi)$.

\noindent  We denote $\hbox{G}(a,b)$, $\hbox{B}(a,b)$, and $\hbox{B}2(a,b, \sigma)$, for $a,b, \sigma >0$, as Gamma, Beta, and Beta type II distributions, and with densities $\frac{\, (bt)^{a} e^{-bt}}{t \Gamma(a) } \, I_{(0,\infty)}(t)$,  $\frac{\Gamma(a+b)}{\Gamma(a)  \, \Gamma(b)} \, 
t^{a-1} \, (1-t)^{b-1} \, I_{(0,1)}(t)$, \hbox{ and }$\frac{\Gamma(a+b) \, \sigma^b}{\Gamma(a)  \, \Gamma(b)} \, 
\frac{t^{a-1}}{(\sigma+t)^{a+b}} \, I_{(0,\infty)}(t)$, respectively.  The latter Beta type II family includes Pareto distributions on $\mathbb{R}_+$ for $a=1$.  

\noindent  The Kummer distribution of type II, denoted $K2(a, b, c, \sigma)$ for parameters $a, c, \sigma >0$ and $b \in \mathbb{R}$, is taken with density
\begin{equation}
\nonumber  \frac{\sigma^{b}}{\Gamma(a) \, \psi(a, 1- b, c)} \, \frac{t^{a-1}}{(t+\sigma)^{a+b}} \, e^{-c t /\sigma} \, I_{(0,\infty)}(t)\,,
\end{equation}
where $\psi$ is the confluent hypergeometric function of type II defined for $\gamma_1, \gamma_3>0$ and $\gamma_2 \in \mathbb{R}$ as:  
$\psi(\gamma_1, \gamma_2, \gamma_3) \, = \, \frac{1}{\Gamma(\gamma_1)} \, \int_{\mathbb{R}_+} t^{\gamma_1 - 1} \, (1+t)^{\gamma_2-\gamma_1-1} \, e^{-\gamma_3 t} \, dt$.  This class of distributions include Gamma distributions with choices $b=-a$.  The class can also be extended to include cases $c=0, b>0$ which correspond to  $\hbox{B}2(a,b, \sigma)$ distributions.
 
\noindent  We will denote McKay's bivariate gamma distribution with parameters $\gamma_1, \gamma_2, \gamma_3 (>0)$ as $Z=(Z_1, Z_2) =^{\!\!\!^d} \hbox{M}(\gamma_1, \gamma_2, \gamma_3)$ with p.d.f. 
\begin{equation}
\nonumber  
\frac{{\gamma_3}^{\gamma_1+\gamma_2}}{\Gamma(\gamma_1) \, \Gamma(\gamma_2)} \, z_1^{\gamma_1-1} \, (z_2-z_1)^{\gamma_2-1} \, e^{-\gamma_3 z_2} \,  \mathbb{I}_{(0,\infty)}(z_1) \, \mathbb{I}_{(z_1,\infty)}(z_2) \,.
\end{equation}
\noindent  The distribution has a long history (e.g., McKay, 1934) and is a benchmark bivariate distribution to model durations that are ordered,   It is easy to verify that the marginals are distributed as $Z_1 =^{\!\!\!^d} \hbox{G}(\gamma_1, \gamma_3)$ and $Z_2 =^{\!\!\!^d} \hbox{G}(\gamma_1+\gamma_2, \gamma_3)$, and that $Z_1$ and $Z_2-Z_1$ are independently distributed, with $Z_2-Z_1 =^{\!\!\!^d} \hbox{G}(\gamma_2, \gamma_3)$.  A generalization will be presented in Section 4.

\section{ Type I mixtures}

\noindent
We begin with Type I mixtures, providing general representations for the posterior distribution $\theta|x$, as well as the predictive distribution $Y|x$, and following up with various examples and observations.   

\begin{theorem}
\label{theoremtypeI}
Let $X,Y|\theta$ be conditionally independent distributed as in (\ref{modeltypeI}) and let $\pi$ have prior density for $\theta$ with respect to $\sigma-$finite measure $\tau$.    
Let $\pi_k$ and $f'_x$ be the densities given by
\begin{equation}
\nonumber
\pi_k(u) \, =\, \frac{\pi(u) \, g_u(k)}{m_{\pi}(k)} \, \hbox{ and }
 f'_x(k) \, = \, \frac{f_k(x) \, m_{\pi}(k)}{\int_{\cal{K}} f_k(x) \, m_{\pi}(k) \, d\mu(k)}\,,\, u \in \Theta, x \in \mathbb{R}^d,
\end{equation}
with $m_{\pi}$ the density given by
\begin{equation}
\nonumber
m_{\pi}(k)= \int_{\Theta} g_{\theta}(k) \, \pi(\theta) \, d\tau(\theta)\,,\, k \in \cal{K}.
\end{equation}
Then,
\begin{enumerate}
\item[ \bf{(a)}]  The posterior distribution of $U =^{\!\!\!^d} \theta|x$ admits the mixture representation: 
\begin{equation}
\label{posttypeI}   U|K' \sim \pi_{K'} \,, \, K' \sim f'_x\,;
\end{equation}
\item[ \bf{(b)}]   The Bayes predictive density of $Y$ admits the representation:
\begin{equation}
\label{predtypeI}  Y|J' \sim q_{J'}\,, J' \sim h'_x,
\end{equation}
with
\begin{equation}
\label{h'}
h'_x(j') \, = \, \int_{\cal{K}} \, q_{\pi}(j'|k') \, f'_{x}(k') \, d\mu(k')\,,\, \hbox{ and }
q_{\pi}(j'|\, k')=\int_{\Theta} h_{\theta}(j') \, \pi_{k'}(\theta) \, d\tau(\theta)\,.
\end{equation}
\end{enumerate}
\end{theorem}
\noindent  {\bf Proof.}
\noindent {\bf (a)} We have indeed
\begin{eqnarray*}
\pi(\theta|x) \, & \propto & \, \int_{\cal{K}}  f_k(x) \, g_{\theta}(k) \, \pi(\theta) \, d\mu(k) \\
\,  & \propto &  \int_{\cal{K}}  f_k(x) \, m_{\pi}(k) \, \pi_k(\theta) \, d\mu(k) \propto \int_{\cal{K}}  f'_k(x) \,\pi_k(\theta)\,  d\mu(k).
\end{eqnarray*}

\noindent {\bf (b)}  The predictive density of $Y$, i.e. the conditional density of $Y|x$, is given by:
\begin{eqnarray*}
q_{\pi}(y|x) \, & = &\, \int_{\Theta}  q(y|\theta) \, \pi(\theta|x) \, d\tau(\theta)  \\
\,  & = & \int_{\Theta} \left(\int_{\cal{K}} q_{j'}(y) \, h_{\theta}(j') \, d\mu(j') \,\right) \,
\left(\int_{\cal{K}} \pi_{k'}(\theta) \, f'_{x}(k') \, d\mu(k') \, \right) \, d\tau(\theta) \\
 & = & \int_{\cal{K}}  q_{j'}(y) \, \left(\int_{\cal{K}} f'_{x}(k') \, \int_{\Theta} \, h_{\theta}(j') \, \pi_{k'}(\theta) \, d\tau(\theta)\ d\mu(k') \right) \, d\mu(j')\,,
\end{eqnarray*}
where we have used (\ref{posttypeI}).  This establishes the result.  \qed

\begin{remark}
\label{interpretation}
The posterior and predictive distribution representations of Theorem \ref{theoremtypeI}
are particularly appealing.   Indeed, observe that posterior distribution (\ref{posttypeI}) mixes the $\pi_{k'}'s$ which correspond to the posterior density of $\theta|X=k'$ as if one had actually observed $K=k'$.  Moreover, the mixing density $f'_{x}$ for $K'$ is a weighted version of the marginal density $m_{\pi}$ for $K$ with weight proportional to $f_k(x)$.

\noindent  The predictive distribution for $Y$ mixes the same densities $q_{J'}$ as the model density for $Y$, with the prior mixing density $h_{\theta}$ replaced by the posterior mixing density $h'_x$.  Furthermore, this mixing density is itself a mixture of the predictive (or conditional) densities $q_{\pi}(\cdot|k')$ of $J$ as if one had observed $K=k'$.
\end{remark}

\noindent  The following examples concern posterior and predictive distribution illustrations of Theorem \ref{theoremtypeI}.  

\begin{example}  
\label{ex1}
In model (\ref{modeltypeI}), consider Poisson mixing $K|\theta =^{\!\!\!^d} \hbox{P}(\theta)$ with a $\hbox{G}(a,b)$ prior for $\theta$.  From this familiar set-up, we obtain $\pi_{k'}$ as a $\hbox{G}(a+k',1+b)$ density and $m_{\pi}$ as a NB$(a, b/(1+b))$ p.m.f.   Following (\ref{posttypeI}), the posterior distribution $\pi(\theta|x)$ is a mixture of the above $\pi_{k}$'s with mixing density on $\mathbb{N}$ given by 
$$f_x'(k) \propto  \frac{(a)_k}{k!} \frac{b^a}{(1+b)^{a+k}} \, f_k(x)\,.  $$  

\noindent Now, consider the cases: {\bf  (i)}  $f_k =^{\!\!\!^d} \hbox{G}(p/2+k,1/2)$, {\bf  (ii)} $f_k =^{\!\!\!^d} \hbox{B}(p/2+k,q/2)$ and {\bf  (iii)} $f_k =^{\!\!\!^d}  \hbox{B}2(p/2+k,q/2)$.  In the context of model (\ref{modeltypeI}), case {\bf (i)} corresponds to a non-central chi-square distribution with $p$ degrees of freedom and non-centrality parameter $2\theta$ 
($\chi^2_p(2\theta)$),  case {\bf  (ii)} to a non-central Beta distribution with shape parameters $p/2$, $q/2$, and non-centrality parameter $2\theta$, and case {\bf  (iii)} to the density of $\frac{p}{q} F$, with $F$ distributed as a non-central $F$ distribution with degrees of freedom $p$, $q$ and non-centrality parameter $2\theta$.  The latter two cases are essentially equivalent though and related by the fact that $\frac{X}{1-X}|\theta$ in  {\bf (ii)} is distributed as $X|\theta$ is in {\bf (iii)}.

\noindent  For {\bf  (i)}, we obtain 
\begin{equation} 
\label{f'(i)}
f_x' \sim \hbox{Hyp}\left(a;\frac{p}{2};\frac{x}{2(1+b)}\right)\,.
\end{equation}  
Observe that the above generalized hypergeometric distribution reduces to a $\hbox{P}(\frac{x}{2(1+b)})$ distribution for $a=p/2$.  The posterior expectation can be evaluated with the help its mixture representation and standard calculations involving the above p.m.f.  One obtains
\begin{eqnarray}
\nonumber \mathbb{E}(\theta|x) \, & = & \mathbb{E}^{K'|x} \{\mathbb{E}(\theta|K',x)  \}
\\
\label{posteriorexpectation}\,  & = & \mathbb{E}_{f'_x}  (\frac{a+K'}{1+b}) \\
\noindent \, & = & \frac{a}{1+b} \left(1 \, + \, \frac{x}{p(1+b)} \, \frac{_1F_1(a+1; p/2+1; \frac{x}{2(1+b)})}{_1F_1(a; p/2; \frac{x}{2(1+b)})}   \right)\,,
\end{eqnarray}
with the case $a=p/2$ simplifying to $\mathbb{E}(\theta|x) \, = \, \frac{p/2}{1+b} \, + \, \frac{x/2}{(1+b)^2}$ as noted by Saxena and Alam (1982).

\noindent For the two other cases, we obtain the mixing power series densities:
\begin{equation}
\label{f'(ii)+f'(iii)}
f_x' \sim \hbox{Hyp}\left(a,\frac{p+q}{2}; \frac{p}{2}; \frac{x}{1+b}\right) \hbox{ for }  {\bf  (ii)}, \hbox{ and }
f_x' \sim \hbox{Hyp}\left(a,\frac{p+q}{2}; \frac{p}{2}; \frac{x}{(1+x)(1+b)}\right) \hbox{ for } {\bf  (iii)}.
\end{equation}

\noindent Observe that the case $a=p/2$ reduces to a NB$(\frac{p+q}{2}, 1 - \frac{x}{1+b})$  distribution for {\bf (ii)} and a NB$(\frac{p+q}{2}, 1 - \frac{x}{(1+x)(1+b)})$  distribution for {\bf (iii)}.  The posterior expectation may be computed from (\ref{posteriorexpectation}) with simplifications occurring for $a=p/2$, yielding for instance in {\bf (ii)} :  $  \mathbb{E}(\theta|x) \, = \, \frac{p(1+b)+qx}{2(1+b-x)}$. Finally, we point out that the above posterior distribution representation applies as well for the improper prior choice $\pi(\theta) = \theta^{a-1} \, \mathbb{I}_{(0,\infty)}(\theta)$ (i.e., $b=0$), with $m_{\pi}$ not a p.m.f., but given by $m_{\pi}(k) \propto \frac{(a)_k}{k!}$.

\end{example}

\begin{example}
\label{expred1}
 Turning now to predictive densities with the same set-up as in Example \ref{ex1}, we consider $Y|\theta$ distributed identically to $X|\theta$ (i.e,  $q_j=f_j$ and $g_{\theta} = h_{\theta}$).  For $Y|\theta =^{\!\!\!^d} \chi^2_p(2\theta)$ (i.e., case {\bf (i)}), Theorem \ref{theoremtypeI} tells us that the Bayes predictive density of $Y$ admits the mixture representation:  
 $Y|J' =^{\!\!\!^d} \hbox{G}(p/2+J',2)$ with $J' \sim h'_x$ given in (\ref{h'}).   The latter admits itself the mixture representation
\begin{equation} 
\label{ex2.3mix}
J'|K' =^{\!\!\!^d} \hbox{NB}(a+K', \frac{1+b}{2+b}),\,\, K' \sim f'_x \hbox{ as in (\ref{f'(i)})}\,,
\end{equation} 
with the former being the Bayes predictive density of $J' =^{\!\!\!^d} \hbox{P}(\theta)$ based on 
$K' =^{\!\!\!^d} \hbox{P}(\theta)$ and prior $\theta =^{\!\!\!^d} \hbox{G}(a,b)$.   Alternatively represented, the mixing p.m.f. for $J'$ may be expressed as
\begin{eqnarray*}
h'_x(j') \,& =& \, \sum_{k'}  \frac{(a+k')_{j'}}{j'!} \, (\frac{1+b}{2+b})^{a+k'} \, (\frac{1}{2+b})^{j'} \, \frac{(a)_{k'}}{k'!\,(\frac{p}{2})_{k'}}  \;  \frac{\left( \frac{x}{2(1+b)}\right)^{k'}}{_1F_1(a; p/2; \frac{x}{2(1+b)})} \\
\,& =& \, \sum_{k'}  \frac{(a+j')_{k'}}{j'!} \, (\frac{1+b}{2+b})^{a} \, (\frac{1}{2+b})^{j'} \, \frac{(a)_{j'}}{k'!\,(\frac{p}{2})_{k'}}  \;  \frac{\left( \frac{x}{2(2+b)}\right)^{k'}}{_1F_1(a; p/2; \frac{x}{2(1+b)})} \\
\, & = & \,  \frac{(a)_{j'}}{j'!} \, (\frac{1+b}{2+b})^{a} \, (\frac{1}{2+b})^{j'} \,
\frac{_1F_1(a+j'; p/2; \frac{x}{2(2+b)})}{_1F_1(a; p/2; \frac{x}{2(1+b)})} \,,
\end{eqnarray*}
which also can be viewed directly as a weighted $\hbox{NB}(a,\frac{1+b}{2+b})$ p.m.f.

\noindent  The non-central Beta and Fisher distributions are similar.  For instance, in the former case with $X|\theta$ and $Y|\theta$ identically distributed with $Y|J =^{\!\!\!^d} \hbox{B}(p/2+J, q/2)$ and $J =^{\!\!\!^d} \hbox{P}(\theta)$, predictive densities associated with $\hbox{G}(a,b)$ priors are also distributed as mixtures
$$  Y|J' =^{\!\!\!^d} \hbox{B}(p/2+J', q/2), \hbox{ with } J'|K' =^{\!\!\!^d} \hbox{NB}(a+K', \frac{1+b}{2+b}),\,\, K' \sim f'_x \hbox{ as in (\ref{f'(ii)+f'(iii)})}\,.$$
For the mixing p.m.f. of $J'$, a development as above yields the expression:

$$  h'_x(j') \, = \,  \frac{(a)_{j'}}{j'!} \, (\frac{1+b}{2+b})^{a} \, (\frac{1}{2+b})^{j'} \,
\frac{_2F_1(a+j',(p+q)/2;  p/2; \frac{x}{2+b})}{_2F_1(a, (p+q)/2 ; p/2; \frac{x}{1+b})} \,.$$
\end{example}

\begin{example}
\label{exampledoublenon-centralF}
A doubly non-central $F$ distribution is a type I mixture (e.g., Bulgren, 1971) with $X \in \mathbb{R_+}$ and $K=(K_1, K_2) \in \mathbb{N}^2$ admitting the representation
\begin{equation}
\label{NCFD}  
X|K =^{\!\!\!^d} \hbox{B}2 \left(\frac{\nu_1}{2}+K_1, \frac{\nu_2}{2}+K_2, \frac{\nu_2}{\nu_1}\right) \hbox{ with } K_1, K_2|\theta =^{\!\!\!\!^{i.d}} \hbox{P}(\frac{\theta_1}{2}) \hbox{ and } \hbox{P}(\frac{\theta_2}{2})\,. 
\end{equation}
Such a distribution arises naturally as a multiple of the ratio of two independent noncentral chi-squared distributions, and reduces to a non-central $F$ for $\theta_2=0$.  Consider now the application of Theorem \ref{theoremtypeI} for the prior 
$\theta_i =^{\!\!\!\!^{i.d}} \, G(a_i,b_i), i=1,2,$ with $a_1,a_2, b_1, b_2>0$, where we have the familiar representations of $\pi_k(\theta_1, \theta_2)$ as the joint density of independent 
$\hbox{G}(a_i+k_i, 1+b_i)\,, i=1,2,$ distributions, and of $m_{\pi}(k_1, k_2)$ as the joint p.m.f.  of independent 
$\hbox{NB}(a_i, \frac{b_i}{1+b_i})\,, i=1,2,$ distributions.
Representation (\ref{posttypeI}) tells us that the posterior distribution of $(\theta_1, \theta_2)$ is a mixture of the $\pi_K$'s with   mixing density  
\begin{eqnarray}
\nonumber 
f'_x(k) \ & \propto &  f_k(x) \, m_{\pi}(k) \\
\nonumber \,  \ & \propto &  \frac{(\frac{\nu_1+\nu_2}{2})_{k_1+k_2} \, (a_1)_{k_1} \, (a_2)_{k_2}}{(\frac{\nu_1}{2})_{k_1} \, (\frac{\nu_2}{2})_{k_2} \, k_1! \, k_2!} \,
\left( \beta_1 \right)^{k_1}  \; \left( \beta_2 \right)^{k_2} \;         I_{\mathbb{N}^2}(k_1,k_2) \\
\label{Appellpmf}
\Rightarrow f'_x(k) & = &  \frac{(\frac{\nu_1+ \nu_2}{2})_{k_1+k_2} \, (a_1)_{k_1} \, (a_2)_{k_2}}{F_2(\frac{\nu_1+\nu_2}{2}; a_1, a_2; \frac{\nu_1}{2}, \frac{\nu_2}{2}; \beta_1, \beta_2)} \;\frac{1}{(\frac{\nu_1}{2})_{k_1} \, (\frac{\nu_2}{2})_{k_2} \, k_1! \, k_2!} \,
\left( \beta_1 \right)^{k_1}  \left( \beta_2 \right)^{k_2} \,         I_{\mathbb{N}^2}(k),
\end{eqnarray}
with $\beta_1 \, = \, \frac{x}{(\frac{\nu_2}{\nu_1} + x)\,(1+b_1)}$, $\beta_2 \, = \, \frac{\nu_2/\nu_1}{(\frac{\nu_2}{\nu_1} + x)\,(1+b_2)}$, and where $F_2$ is the Appell function of the second kind given by
$$ F_2(\gamma_1; \gamma_2,  \gamma_3; \gamma_4, \gamma_5; w,z) \, = \, 
 \sum_{m=0}^{\infty} \sum_{n=0}^{\infty} \, \frac{(\gamma_1)_{m+n} \, (\gamma_2)_m \, (\gamma_3)_n}{m! \, n! \, (\gamma_4)_m \, (\gamma_5)_n} \, w^m \, z^n\,.$$
 
\noindent  The bivariate p.m.f. $f'_x$ in (\ref{Appellpmf}) and how it is arisen here are of interest.  It is a bivariate power series p.m.f. generated by the coefficients of the Appell $F_2$ function and will be bona fide p.m.f. for $\beta_1, \beta_2>0$ and $\beta_1+\beta_2 < 1$.   Appell's $F_1$ function appears in a similar way, again in a Bayesian framework, as a bivariate discrete distribution called {\it Bailey} by Laurent (2012) (see also Jones \& Marchand, 2019 for another derivation).  For the particular case $a_i=\nu_i/2, i=1,2$, the p.m.f. in (\ref{Appellpmf}) simplifies and the corresponding random pair $(K_1', K_2')$, admits the stochastic representation:  $K_1' \sim \hbox{NB}(a_1+a_2, 1-p_1)$ and $K_2'|K_1' \sim \hbox{NB}(a_1+a_2+K_1', p_2)$ with $p_1= \beta_1/(1-\beta_2) $ and $p_2=1-\beta_2$.

\normalsize
\noindent  Turning to the predictive density, consider $Y|\theta$ distributed as $X|\theta$ and the same prior on $\theta$ as above.   Similarly to Example \ref{expred1}, which shares the same Poisson mixing and gamma prior, we obtain from Theorem \ref{theoremtypeI} and (\ref{predtypeI}), that the predictive density for $Y$ admits the representation
\begin{equation}
\nonumber Y|J' =^{\!\!\!^d} \hbox{B}2 \left(\frac{\nu_1}{2}+J'_1, \frac{\nu_2}{2}+J'_2, \frac{\nu_2}{\nu_1}\right), \hbox{ with } J'_i|K' =^{\!\!\!^{i.d}} \hbox{NB}(a+K'_i, \frac{1+b_i}{2+b_i}),\,\, K' \sim f'_x \hbox{ as in (\ref{Appellpmf})}\,.
\end{equation}
\noindent  We point out that, if the distribution of $Y|\theta$ is of the same type as that of $X|\theta$ (in \ref{NCFD}) with associated degrees of freedom $\nu_1'$ and $\nu_2'$ instead of $\nu_1$ and $\nu_2$, the only change in the previous expression is to replace the $\nu_i$'s by the $\nu_i'$'s for the distribution of $Y|J'$.  When applicable, a similar observation applies to the other examples of this section.  

\end{example}

\begin{example}
\label{exgammaminxtures}
Consider univariate Gamma mixtures with $X|K =^{\!\!\!^d} G(\alpha, K)$ in model (\ref{modeltypeI}), with the mixing distribution $K|\theta =^{\!\!\!^d} G(a,\theta)$ and prior 
$\theta =^{\!\!\!^d} G(c,b)$, with known $\alpha, a, b, c$.  Theorem \ref{theoremtypeI} applies and tells us that the posterior distribution $\theta|x$ is a mixture of $\pi_{K'} \sim G(a+c, K'+b)$ making use of a familiar posterior distribution for gamma models with gamma priors.  In evaluating the mixing density $f'_x$ of $K'$ given in (\ref{posttypeI}), it is easy to verify that $m_{\pi} \sim \hbox{B}2(a,c,b)$ and one thus obtains
$$ f'_x(k) \propto  m_{\pi}(k) \, f_{k}(x)  \propto \frac{k^{a+\alpha-1}}{(k+b)^{a+c}} \, e^{-xk} \; \mathbb{I}_{(0,\infty)}(x)\,,$$
which corresponds to a Kummer $\hbox{K}2(a+\alpha, c-\alpha, bx, b)$ distribution as defined in Section 2. 

\noindent  Now consider the Bayesian predictive density for the Gamma mixture $Y|J =^{\!\!\!^d} \hbox{G}(\alpha', J)$ and $J|\theta =^{\!\!\!^d} \hbox{G}(a',\theta)$, which includes the particular case of identically distributed $Y|\theta$ and 
$X|\theta$ for $\alpha'=\alpha$ and  $a'=a$.  Observe that the density $q_{\pi}(\cdot|k')$ is that of $J|K=k'$ for the set-up $J|\theta =^{\!\!\!^d} \hbox{G}(a',\theta)$, 
$K|\theta =^{\!\!\!^d} \hbox{G}(a,\theta)$ independently distributed, and with $\theta =^{\!\!\!^d} \hbox{G}(c,b)$.  A calculation (e.g., Aitchison \& Dunsmore, 1975) yields $ q_{\pi}(\cdot|k') =^{\!\!\!^d} \hbox{B2}(a', a+c, b+k')$.   With the above, it follows from Theorem \ref{theoremtypeI} that the predictive distribution for $Y$ admits the representation:
\begin{equation}
\nonumber Y|J' =^{\!\!\!^d} \hbox{G}(\alpha', J')\,, \hbox{ with }J'|K' =^{\!\!\!^d} \hbox{B2}(a', a+c, b+K') \, \hbox{ and } \, 
K' =^{\!\!\!^d} \hbox{K}2(a+\alpha, c-\alpha, bx, b)\,.
\end{equation}
An alternative representation comes from simply evaluating the marginal density $h'_x$ of $J'$.  A calculation gives:
\begin{equation}
\nonumber h'_x(j') \, = \,  \frac{\Gamma(a'+a+c)}{\Gamma(a') \, \Gamma(a+c) } \,
\frac{b^{c-\alpha} \, j'^{(a'-1)}}{(j'+b)^{a'+c-\alpha}} \, \frac{\psi\left( a+\alpha, 1+\alpha-a'-c, jx+bx \right)}{\psi\left(a+\alpha, 1+\alpha-c, bx\right)} \; 
 \mathbb{I}_{(0,\infty)}(j')\,,
\end{equation}
with $\psi$ the confluent hypergeometric function of type II as defined in Section 2. 
\end{example}

\noindent  The final application of Theorem \ref{theoremtypeI} concerns the coefficient of determination in a standard multiple regression context.

\begin{example}    Consider a coefficient of determination $X=R_1^2$, or square of a multiple correlation coefficient, that arises for a sample of size $n_1$ from $Z=(Z_1, \ldots, Z_m)^{\top}$, with $n_1>m>1$, and the regression of $Z_1$ based on $Z_2, \ldots, Z_m$.     For more details on the underlying distributional theory, see for instance Muirhead (1982).   It is well known that the distribution of $X$ is a Type I mixture (\ref{modeltypeI}) with
\begin{equation}
\nonumber   X|K =^{\!\!\!^d} \hbox{B}\left(\frac{m-1}{2}+K, \frac{n_1-m+1}{2}\right)\,, \hbox{ with }  K|\theta =^{\!\!\!^d} \hbox{NB}\left(\frac{n_1-1}{2}, 1-\theta\right)\,,
\end{equation}
with $\theta=\rho^2$ is the theoretical squared multiple correlation coefficient.  As in Marchand (2001), a convenient prior on $\theta$ is a Beta prior and it leads along with the negative binomial distributed $K|\theta$ to a conjugate posterior.   Specifically, we obtain for $\theta =^{\!\!\!^d} \hbox{B}(a,b)$ the posterior density $\pi_K \sim \hbox{B}\left(a+K, b+(n_1-1)/2\right)$, as well as the marginal p.m.f.
\begin{eqnarray*}
m_{\pi}(k) \, &  = &\, \int_{(0,1)} \, \frac{\Gamma(a+b)}{\Gamma(a) \Gamma(b)} \,\theta^{a-1} (1-\theta)^{b-1}
\,   \frac{\Gamma(\frac{n_1-1}{2}+k )}{k! \, \Gamma(\frac{n_1-1}{2})} \, \theta^k \, 
(1-\theta)^{(n_1-1)/2}  \, d\theta \\
\, & = & \frac{(\frac{n_1-1}{2})_k \, (a)_k}{k! \, (a+b+\frac{n_1-1}{2})_k} \, \frac{(b)_{(n_1-1)/2}}{(a+b)_{(n_1-1)/2}}  \,, \hbox{ for } k \in \mathbb{N}\,.
 \end{eqnarray*}

 Theorem \ref{theoremtypeI} tells us that the posterior distribution $\theta|x$ is a mixture of the $\pi_k$'s with mixing 
\begin{eqnarray}
\nonumber f'_x(k) \, & \propto & \, f_k(x) \, m_{\pi}(k) \\
\nonumber\,   & \propto &    \frac{(\frac{n_1-1}{2})_k}{(\frac{m-1}{2})_k} \,x^k \, \frac{(\frac{n_1-1}{2})_k \, (a)_k}{k! \, (a+b+\frac{n_1-1}{2})_k} \\
\label{f'R^2} \,  
\Rightarrow  f'_x \, & \sim & \, \hbox{Hyp}\left(\frac{n_1-1}{2}, \frac{n_1-1}{2}, a; \frac{m-1}{2}, a+b+ \frac{n_1-1}{2}; x  \right) \,.
\end{eqnarray} 
The result, which we have derived from the general context of Theorem \ref{theoremtypeI} was obtained by Marchand (2001) in this specific set-up.  In doing so, he defined such Beta mixtures as HyperBeta and also provided several graphs of prior-posterior densities for varying prior parameters $a,b$, sample size $n$, and observed values of $X$.  \\

\noindent For the predictive density of a future $Y=R^2_2$ distributed as $X$ but allowing for a possibly different sample size $n_2$,  expression (\ref{predtypeI}) tells us that such a predictive density admits the mixture representation:
\begin{equation}
\label{predR^2}
Y|J' =^{\!\!\!^d} \hbox{B}(\frac{m-1}{2}+J', \frac{n_2-m+1}{2})\,, J'|k' \sim q_{\pi}(\cdot|k'), \, \hbox{ and } K' \sim f'_x \hbox{ as in } (\ref{f'R^2})\,, with
\end{equation}
$q_{\pi}(\cdot|k')$ the predictive density for $J'|\theta =^{\!\!\!^d} \hbox{NB}(\frac{n_2-1}{2}, 1-\theta)$ based on $K'|\theta =^{\!\!\!^d} \hbox{NB}(\frac{n_1-1}{2}, 1-\theta)$ and prior $\theta =^{\!\!\!^d} \hbox{B}(a,b)$.  An evaluation of (\ref{h'}) yields  
\begin{eqnarray*}
q_{\pi}(j'|k')\, & = &\, \int_0^1  \frac{(\frac{n_2-1}{2})_{j'}}{j'!} \, \theta^{j'} \, (1-\theta)^{(n_2-1)/2} \, \frac{\Gamma(a+k'+b+\frac{n_1-1}{2})}{\Gamma(a+k') \, \Gamma(b + \frac{n_1-1}{2})} \, \theta^{a+k'-1} (1-\theta)^{b+(n_1-1)/2-1}   d\theta \\ 
\, & = & \, \frac{\Gamma(b+(n_1+n_2)/2 -1) \, \Gamma(a+k'+b+\frac{n_1-1}{2}) }{\Gamma(a+k'+b+(n_1+n_2)/2-1) \, \Gamma(b+\frac{n_1-1}{2})} \; \frac{(\frac{n_2-1}{2})_{j'} \, (a+k')_{j'}}{j'! \, (a+k'+b+(n_1+n_2)/2 -1)_{j'}} \,,  \hbox{ for } j' \in \mathbb{N}.
\end{eqnarray*}
It is worth noting that the above belongs to the class of p.m.f.'s of the form $p(n)\, = \frac{1}{_2F_1(a,b;c,1) }\, \frac{(a)_n \, (b)_n}{(c)_n \, n!} \, \mathbb{I}_{\mathbb{N}}(n)$ with  $a,b,c>0$, $c>a+b$, and as ${}_2F_1(a, b; c; 1)\, = \, \frac{\Gamma(c)\,\Gamma(c-a-b)\,}{\Gamma(c-a)\,\Gamma(c-b)\,}$.

\end{example}

\noindent As a complement to the above, we briefly illustrate the usefulness of the predictive density (\ref{predR^2}) for visualizing how concordant or discordant the sample coefficients of determination $X=R_1^2$ and $Y=R_2^2$ can be when arising from studies set under the same conditions.  Such a comparison also relates to the degree to which a second study replicates or not the results of an initial study when focusing on the $R^2$ summaries.   Consider Figure 1 where both posterior densities for $\rho^2$ and predictive densities for $Y=R_2^2$ are drawn for $X=R_1^2=0.25, 0.50$, a $\hbox{B}(3,5)$ prior, and sample sizes $n_1=27, n_2=31$.
The variability in the Bayesian predictive distributions is indeed quite significant, stemming of course from posterior uncertainty as well as model uncertainty.    Further illustrating and for further interpretation of Figure 1, we record the expectations and standard deviations under the predictive densities:  $\mathbb{E}(R_2^2|R_1^2=0.25) \, \approx 0.323 \,,$ 
$\sigma(R_2^2|R_1^2=0.25) \, \approx 0.159 \,,$  $\mathbb{E}(R_2^2|R_1^2=0.50) \, \approx 0.431 \,, $  and $\sigma(R_2^2|R_1^2=0.50)\,, \approx 0.170$,
as well as the selected probabilities  $\mathbb{P}(R_2^2>0.40|R_1^2=0.25) \, \approx 0.310 \, $ and 
$\mathbb{P}(R_2^2<0.10|R_1^2=0.25) \, \approx 0.087 \, $.

 \begin{figure}[ht]
 \label{graphe1}
 \centering
\includegraphics[width=0.79\textwidth]{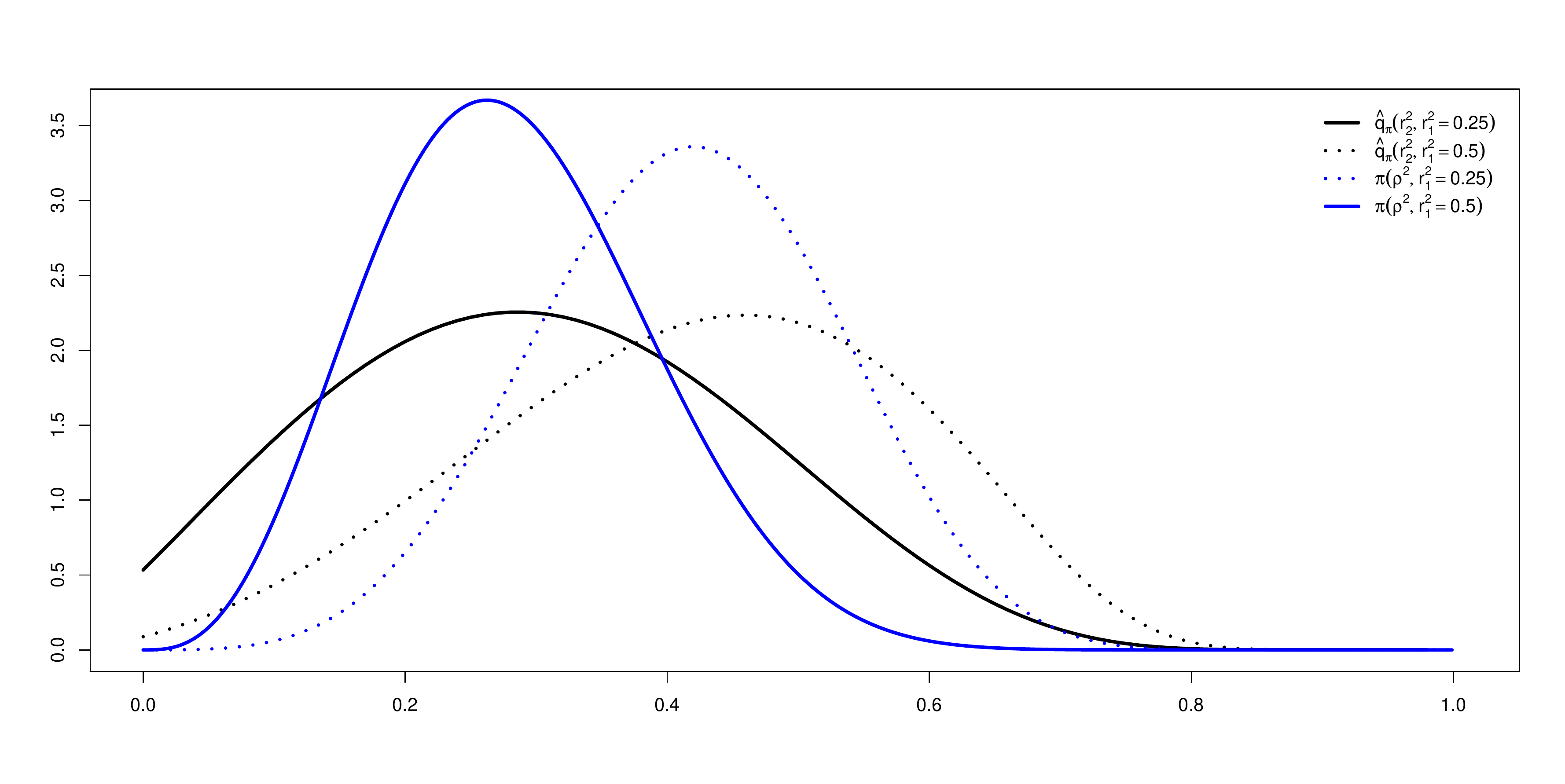}
  \caption{Bayesian predictive densities $q_{\pi}$ for a future $R_2^2$ ($n_2=31)$ based on two values of $R_1^2$ ($n_1=27$), and a  $\hbox{B}(3,5)$ prior. } \label{graphe1}
 \end{figure}

\section{Type II mixtures}

\noindent
We consider now Type II mixtures, provide general representations for the posterior distribution $\theta|x$, as well as the predictive distribution $Y|x$, and continue with various examples and observations.

\begin{theorem}
\label{theoremtypeII}
Let $X,Y|\theta$ be conditionally independent distributed as in (\ref{modeltypeII}) and let $\pi$ have prior density for $\theta$ with respect to $\sigma-$finite measure $\tau$.    
Let $\pi_{k,x}$ and $g_{\pi,x}$ be respectively the conditional distribution densities of $\theta|x,k$ and $k|x$ given by 
\begin{equation}
\nonumber
\pi_{k,x}(\theta) \, =\, \frac{\pi(\theta) \, f_{\theta,k}(x)}{m_{\pi,k}(x)} \, \hbox{ and }
 g_{\pi,x}(k) \, = \, \frac{g(k) \, m_{\pi,k}(x)\,}{\int_{\cal{K}} g(k) \, m_{\pi,k}(x) \, d\mu(k)},\, \theta \in \Theta, x \in \mathbb{R}^d, k \in \cal{K}\,,
\end{equation}
with $m_{\pi,k}$ given by
\begin{equation}
\nonumber
m_{\pi,k}(x)= \int_{\Theta} f_{\theta,k}(x) \, \pi(\theta) \, d\tau(\theta)\,.
\end{equation}
Then,
\begin{enumerate}
\item[ \bf{(a)}]  The posterior distribution $U =^{\!\!\!^d} \theta|x$ admits the mixture representation: 
\begin{equation}
\label{posttypeII}   U|K' \sim \pi_{K',x} \,, \, K' \sim g_{\pi,x}\,;
\end{equation}

\item[ \bf{(b)}]   The Bayes predictive density of $Y$ admits the representation:
\begin{equation}
\label{preddensitytypeII}  Y|J',K' \sim q_{\pi}(\cdot|J', K')\,, J' \sim h\,,\,K' \sim g_{\pi,x} 
{\hbox{ independent}}\,,
\end{equation}
with
\begin{equation}
\label{qpi}
 q_{\pi}(y|\, j', k')=\int_{\Theta} q_{\theta,j'}(y) \, \pi_{k',x}(\theta) \, d\tau(\theta)\,.
\end{equation}
An alternative equivalent representation is: 
\begin{equation}
\nonumber
\label{predtypeII}  Y|J' \sim q'_{x,J'}\,, J' \sim h,
\end{equation}
with
\begin{equation}
\nonumber
\label{q'}
q'_{x, j'}(y) \, = \, \int_{\cal{K}} \, q_{\pi}(y|j',k') \, g_{\pi,x}(k') \, d\mu(k')\,. 
\end{equation}
\end{enumerate}
\end{theorem}
\noindent  {\bf Proof.}
{\bf (a)} 
We have indeed
\begin{equation}
\nonumber
\pi(\theta|x) \, = \, \int_{\cal{K}}  \, \pi(\theta|x,k) \, g_{\pi,x}(k)  \, d\mu(k) \,,
\end{equation}
with
\begin{equation}
\nonumber
\pi(\theta|x,k) \, = \, \frac{\pi(\theta) \, f_{\theta,k}(x) \, g(k)}{\int_{\Theta} \pi(\theta) \, f_{\theta,k}(x) \, g(k) \, d\tau(\theta)}\, \, = \, \pi_{k,x}(\theta)\,,
\end{equation}

\begin{equation}
\nonumber
\hbox{ and }  \, g_{\pi,x}(k) \, \propto \, \int_{\Theta} f_{\theta,k}(x) \, g(k) \, \pi(\theta) \, d\tau(\theta)\, = \, g(k) \, m_{\pi,k}(x)\,, 
\end{equation}
which establishes (\ref{posttypeII}).

\noindent {\bf (b)}  The Bayes predictive density of $Y$, i.e. the conditional density of $Y|x$, is given by:
\begin{eqnarray*}
q_{\pi}(y|x) \, & = &\, \int_{\Theta} q(y|\theta) \, \pi(\theta|x) \, d\tau(\theta) \\
\,  & = & \,  \int_{\Theta} \left( \int_{\cal{K}} \, q_{\theta,j'}(y) \, h(j') d\mu(j')  \, \right) \, \pi(\theta|x) \, d\tau(\theta) \\
\,  & = & \,  \int_{\Theta} \left( \int_{\cal{K}} \, q_{\theta,j'}(y) \, h(j') d\mu(j') \, \int_{\cal{K}} \pi_{k',x}(\theta) \, g_{\pi,x}(k') \,d\mu(k') \right) \, d\tau(\theta) \\
\, & = & \,  \int_{\cal{K}}  \int_{\cal{K}}  \left( \int_{\Theta} q_{\theta,j'}(y) \, \pi_{k',x}(\theta) \, d\tau(\theta) \right) \,  g_{\pi,x}(k') \, h(j') \, d\mu(k') \, d\mu(j')  \\
\, & = & \,   \int_{\cal{K}}  \int_{\cal{K}}  q_{\pi}(y|j',k') \,  g_{\pi,x}(k') \, h(j') \, d\mu(k') \, d\mu(j')\,,
\end{eqnarray*}
where we have used (\ref{posttypeII}).  This establishes the result.  \qed

\begin{remark}
\label{remarktypeII}  Again, the posterior and predictive distributions have noteworthy representations.  Namely, while the model distribution for $X|\theta$ mixes the $f_{\theta,k}$'s with mixing marginal density $g$ of $K$, the posterior distribution (\ref{posttypeII}) mixes the $\pi_{k',x}$'s, which represents for fixed $k',x$ the posterior density of $\theta$ as if one had actually observed as well $K=k'$, with mixing density $g_{\pi,x}$ which is the marginal posterior distribution of $K$.  Similarly, the density $m_{\pi,k'}$ is the marginal density of $X$ as if $K=k'$ had been observed.

\noindent  In the second representation of the predictive distribution for $Y$, we have a mixture with the same mixing density $h$ as of the model distribution of $Y$.   The distributions mixed are themselves mixtures of the $q_{\pi}(\cdot|j',k')$ densities, which represent Bayes predictive densities for $Y$ as if one had actually observed $K=k'$ and $J=j'$.  With the first representation, we have bivariate mixing of the densities $q_{\pi}(\cdot|J',K')$ with independently distributed components $J'$ and $K'$.

\end{remark}

\begin{example}
\label{examplemixturesnomals}
As a first illustration of the above Theorem, consider variance mixture of normal distributions with summaries $X|\theta,K =^{\!\!\!^d} N_d(\theta, K I _d)$ with $K \sim g$ on $\mathbb{R}_+$, as well as $Y|\theta,J =^{\!\!\!^d} N_d(\theta, J I _d)$ with $J \sim h$ on $\mathbb{R}_+$.  Such a class includes not only finite mixtures with the mixing variables $K$ and $J$ discrete with finite support, but also several familiar distributions such as multivariate Student, Logistic and Laplace, among others.  

\noindent  Now, consider the prior $\theta =^{\!\!\!^d} N_d(\mu, \tau^2 I_d)$ with known $\mu, \tau^2$.  In applying Theorem \ref{theoremtypeII} for the posterior distribution of $\theta$, it is well known that 
\begin{eqnarray}
\label{pi_x,k(normal)}
\nonumber \pi_{k,x} & \sim & N_d\left( x - \frac{k}{\tau^2+k} (x-\mu), \frac{\tau^2 k}{\tau^2+k} I_d\right)\,, \\
\nonumber \hbox{ and  } \; m_{\pi,k} & \sim & N_d\left( \mu, (\tau^2+k) \,  I_d\right).
\end{eqnarray}
Theorem \ref{theoremtypeII} thus implies that $U =^{\!\!\!^d} \theta|x$ admits the representation
\begin{equation}
\label{gpixnormalmix}
U|K' \sim \pi_{K',x} \hbox{ with } K' \sim g_{\pi,x}(k) \propto \frac{g(k)}{(k+\tau^2)^{d/2}} \,\, e^{- \frac{\|x-\mu\|^2}{2(k+\tau^2)}}\,.
\end{equation}
Alternatively, we can express the above in terms of $Z=\frac{\theta-x}{\tau}$, $R = \frac{K'}{\tau^2+K'}$, and $\beta \, = \, \frac{\mu-x}{\tau}$, yielding the mixture representation 
\begin{equation}
\nonumber
Z|R =^{\!\!\!^d} N_d \left(\beta R, R I_d \right), \, \hbox{ with } R \sim g_{R,x}(r) \propto   g\left(\frac{\tau^2 \, r}{1-r}\right) \, (1-r)^{d/2-2} \,\, e^{- \beta^2 (1-r)/2}\, \mathbb{I}_{(0,1)}(r)\,,
\end{equation}
$g_{R,x}$ being the density of $R$. 
The above posterior distribution for $Z$ is an example of a variance-mean mixture of normal distributions popularized by Barndorff-Nielson et al. (1982), and much studied since in terms of theoretical and practical issues.   For the improper prior $\pi(\theta)=1$, it is not too difficult to evaluate directly the posterior, but Theorem \ref{theoremtypeII} also applies.  Indeed in this case, we have $\pi_{k,x} \sim N_d(x, k I_d)\,, m_{\pi,k}(x)=1$, and the posterior distribution $U =^{\!\!\!^d} \theta|x$ admits the representation of $U|K' =^{\!\!\!^d} N_d(x, K I_d)$ with $K' \sim g$, which is a variance (only) mixture of normal distributions.

\noindent  For the predictive distribution of $Y$, one can verify the Bayesian normal model predictive density calculation
\begin{equation}
\nonumber  q_{\pi}(\cdot|J',K') \sim N_d \left(x - \frac{K'}{\tau^2+K'}(x-\mu)\,,\, (\frac{\tau^2 \, K'}{\tau^2+K'}) \, + J') I_d  \right)\,.
\end{equation}
With this, it follows from Theorem \ref{theoremtypeII} that the Bayes predictive density of $Y$ is a mixture of normal densities $q_{\pi}(\cdot|J',K')$ with $J' \sim h$ and $K' \sim g_{\pi,x}$ given in (\ref{gpixnormalmix}).  Observe that the result applies if either $X$ or $Y$ is normally distributed by taking the corresponding mixing variable $K$ or $J$ to be degenerate.  Finally, for the improper prior $\pi(\theta)=1$, the above becomes
$q_{\pi}(\cdot|J',K') \sim N_d(x, (J'+K') I_d)$ with $J' \sim h$ and $K' \sim g$, a result given by Kubokawa et al. (2015).   

\end{example}

\noindent  The next example deals with a mixture of exponential distributions (denoted ${\cal{E}}$) model which has garnered much over the years much interest in reliability.
\begin{example}
\label{lomax}
We now consider here the following mixture of exponential distributions model:
\begin{equation}
\label{modellomax}   
Z_1, \ldots, Z_n|K \,=^{\!\!\!\!^{i.d}}\,  {\cal{E}}(\mu, \frac{K}{\sigma})\,,\, K \sim g \,, 
\end{equation}
with $n \geq 2$, and with unknown location and scale parameters $\mu$ and $\sigma$.  Here, the densities being mixed are exponential distributions such that $X_i-\mu|K \,=^{\!\!\!\!^{d}} \hbox{G}(1, \frac{K}{\sigma})$.  The above is a type II mixture and we set $\theta=(\mu,\rho)$ with $\rho=1/\sigma$.   Quite generally, model (\ref{modellomax}) is representative of $n$ system components with a common environment, the case $n=2$ having been put forth by Lindley \& Singpurwalla (1986).  A notable choice of $g$ occurs for $K \,=^{\!\!\!\!^{d}} \hbox{G}(\alpha,1)$ with known $\alpha$, and leads to the multivariate Lomax distribution (e.g., Nayak, 1987) which has generated much interest over the years.  

\noindent  It is reasonable to consider the summary $(S,R)$ with $S=\sum_{i=1}^n (Z_i - Z_{(1)}) $, $R=Z_{(1)} \, = \, \min_i Z_i$ (e.g., Petropoulos, 2017).  The pair $(S,R)$ is, namely, a jointly sufficient and complete statistic for degenerate $K$.   Conditional on $K$, $S$ and $R$ are independently distributed as:
\begin{equation}
\label{complete+sufficient}   S|K \,=^{\!\!\!\!^{d}} \hbox{G}(n-1, \rho K) \, \hbox{ and } \, R|K \,=^{\!\!\!\!^{d}} {\cal{E}}(\mu, nk \rho)\,.
\end{equation}
We thus have a Type II mixture for $X=(S,R)$ with
\begin{equation}
\label{densitymodellomax}
f_{\theta, k}(x) \, = \, \frac{n (\rho k)^n}{\Gamma(n-1)} \, s^{n-2} \, e^{-\rho k \left\lbrace s+ n(r-\mu)\right\rbrace} \,\, \mathbb{I}_{\mathbb{R}_+}(s) \, \mathbb{I}_{\mathbb{R}_+}(r-\mu)\,. 
\end{equation}
Consider now a gamma prior for $\rho$ and a flat prior for $\mu$ so that
\begin{equation}
\label{priorlomax}
\pi(\rho, \mu) \, \propto \, \rho^{a-1} \, e^{-b \rho} \, \mathbb{I}_{\mathbb{R}_+}(\rho) \; \mathbb{I}_{\mathbb{R}}(\mu)\,, a,b>0\,.
\end{equation}
In fact, the posterior distribution exists as well for $b=0$, and $a \in (2-n, 0]$ for $n>2$.  Evaluating the density $\pi_{k,x}$ of Theorem \ref{theoremtypeII}. we obtain
\begin{equation}
\pi_{k,x}(\rho, \mu) \,\propto \, \{\rho^{a+n-2} \, e^{-\rho(ks+b)} \}\, \{\rho k n \, e^{-\rho k n  (r-\mu)}\,\} \; \mathbb{I}_{\mathbb{R}_+}(\rho) \, \mathbb{I}_{\mathbb{R}_+}(r-\mu)\,,
\end{equation}
which implies the posterior distribution representation:
\begin{equation}
\label{postlomax}  \rho|K,x \,=^{\!\!\!\!^{d}} \hbox{G}(a+n-1, Ks+b)\,,\, \omega|\rho, K, x \,=^{\!\!\!\!^{d}} {\cal{E}}(0, \rho K n)\,,
\end{equation}
where we have set $\omega=r-\mu$.   Now, evaluating the marginal density\footnote{As a function of $x$, it is not a true density here as we have used an improper prior.} $m_{\pi,k}$ of $X=(S,R)$, we have from (\ref{densitymodellomax}) and (\ref{priorlomax}):
\begin{eqnarray*}
m_{\pi,k}(x) & = & \int_{\Theta} \pi(\theta) \, f_{\theta,k}(x) \, d\tau(\theta) \\
\, & \propto & \, k^n \, \int_{\mathbb{R}_+} \, \rho^{a+n -1} \, e^{-\rho(b+ks)} \,  \{\int_{-\infty}^r  \, e^{-\rho nk(r-\mu)}d\mu\} \, d\rho \\
\, &   \propto &  \frac{k^{n-1}}{(\frac{b}{s} + k)^{a+n-1}}\,,
\end{eqnarray*}
as a function of $k \in \cal{K}$.   Finally, putting the above elements together, we have for $X|K,\theta$ as in (\ref{complete+sufficient}) with $K  \sim g$ and prior density (\ref{priorlomax}), the following representation for the posterior distribution of $(\rho, \omega=r-\mu)$:
\begin{equation}
\nonumber \rho|K',x \,=^{\!\!\!\!^{d}} \hbox{G}(a+n-1, K's+b)\,,\, \omega|\rho, K', x \,=^{\!\!\!\!^{d}} {\cal{E}}(0, \rho K' n)\, \hbox{ with } K' \sim g_{\pi,x}(k) \propto\, \frac{g(k) \, k^{n-1}}{(\frac{b}{s} + k)^{a+n-1}}\,.  
\end{equation} 
From this, we obtain for $g \,=^{\!\!\!\!^{d}} \hbox{G}(\alpha, \beta)$ the Kummer density 
\begin{equation} g_{\pi,x} \,=^{\!\!\!\!^{d}} \hbox{K2}(\alpha+n-1, a+1-\alpha, \frac{b \beta}{s}, \frac{b}{s})\,,
\end{equation}
with $\beta=1$ corresponding to the multivariate Lomax case.  Alternatively, for an inverse Gamma density $g$ with density proportional to $k^{-\alpha-1} \, e^{-\beta/k} \, \mathbb{I}_{(0,\infty)}(k)$, we obtain that the mixing density of $K'$ is such that $\frac{1}{K'} \,=^{\!\!\!\!^{d}} \hbox{K}2(a+\alpha+2, n-\alpha-3, \beta s /b, s/b)$.  \\
 
\noindent   Now, for deriving a predictive density via Theorem \ref{theoremtypeII}, based on $X=(S,R)$ distributed as above in (\ref{complete+sufficient}) with $n=n_1$, consider predicting $Y=(S',R')$ such that:
\begin{equation}
S'|J \,=^{\!\!\!\!^{d}} \hbox{G}(n_2-1, \rho J) \,\hbox{ independent of } \, R'|J \,=^{\!\!\!\!^{d}} {\cal{E}}(\mu, n_2 \rho J)\,, \hbox{ with }  J \sim h\,.
\end{equation}
These statistics arise naturally, just as the case for $X=(S,R)$, as jointly complete and sufficient conditional for the observables $W_1, \ldots, W_{n_2}|J \,=^{\!\!\!\!^{i.d}}  {\cal{E}}(\mu, J \rho)\,$ with  $S'=\sum_{i=1}^{n_2} (W_i - W_{(1)})$ and $R'= W_{(1)}$.
Following Theorem \ref{theoremtypeII}, based on a prior density $\pi$ for $\theta=(\rho, \mu)$, we have that the Bayes predictive density of $Y$ admits representation (\ref{preddensitytypeII}), which is $Y|J',K' \sim q_{\pi}(\cdot|J',K')$ with $J' \sim h$ and $K' \sim g_{\pi,x}$.   There thus remains and we now seek to evaluate
\begin{equation}
\label{qpilomaxj'k'}
q_{\pi}(s',r'|j',k') \, = \, \int_{\mathbb{R}^2}  q_{\theta,j'}(s',r') \, \pi_{k',x}(\rho, \mu) \, d\rho \, d\mu \,,
\end{equation}
with
\begin{equation}
\nonumber  q_{\theta,j'}(s',r') \, = \,  \frac{{s'}^{(n_2-2)} \, e^{-\rho j' s'}}{\Gamma(n_2-1)} \, (\rho j')^{n_2} \, n_2 \, e^{-n_2 j' \rho (r'-\mu)} \, \mathbb{I}_{\mathbb{R}_+}(s') \, \mathbb{I}_{(\mu, \infty)}(r') 
\end{equation}
and
\begin{equation}
\nonumber  \pi_{k',x}(\rho, \mu) \, = \,  \frac{\rho^{a+n_1-2} \, e^{-\rho(k's+b)}}{\Gamma(a+n_1-1)} \,(k's+b)^{a+n_1-1} \, (\rho n_1 k') \, e^{-\rho k' n_1  (r-\mu)}\,   \mathbb{I}_{\mathbb{R}_+}(\rho) \, \mathbb{I}_{(-\infty, r)}(\mu)\,.
\end{equation}
\noindent After some calculations, setting $A_0 \, = \frac{n_1 n_2 k' j'^{n_2} \, {s'}^{(n_2-2)} \, \,(k's+b)^{a+n_1-1}}{\Gamma(n_2-1) \, \Gamma(a+n_1-1)} \, $ and $B_0 \, = \, j's' + k's + b + n_2 j' r' + n_1 k' r$,  one obtains
\begin{eqnarray*}
q_{\pi}(s',r'|j',k') \,& = &\, A_0  \int_0^{\infty} \rho^{a+n_1+n_2-1} \, e^{-\rho B_0} \int_{-\infty}^{r \wedge r'} \, e^{\rho \mu (n_2 j' + n_1 k')}  d\rho \, d\mu \\
\, & =  & \,  \frac{A_0 \, }{(n_2 j' + n_1 k')} \; \frac{\Gamma(a+n_1+n_2-1)}{\left\lbrace  j's' \, + \, k' s \, + \,  b \,+ \, c(r') \, \right\rbrace^{a+n_1+n_2-1}} \,\, \mathbb{I}_{\mathbb{R}_+}(s') \,\, \mathbb{I}_{\mathbb{R}}(r')\,,
\end{eqnarray*}
with $c(r') \, = \, |r'-r| \, \left\lbrace n_2 j' \, I_{(r, \infty)}(r') \, + n_1 k' \,   I_{(- \infty, r)}(r') \, \right\rbrace$.  This provides an explicit expression for the sought-after (\ref{qpilomaxj'k'}), but a more appealing decomposition of the distribution of $(S',R')$ under $q_{\pi}(\cdot|j',k')$, easily derived from the above, is:
\begin{eqnarray}
\label{Betalomax}S'|R',J',K' & \,=^{\!\!\!\!^{d}} &  \hbox{B2}(n_2-1, a+n_1, \frac{K's+b+c(r')}{J'}) \\
\nonumber R'|J',K' & \sim &   p  \,  \frac{A_1^{a+n_1-1} (a+n_1-1)}{(A_1+|r'-r|)^{a+n_1}} \, I_{(r, \infty)}(r') \; + \;  (1-p)  \,  \frac{A_2^{a+n_1-1} (a+n_1-1)}{(A_2+|r'-r|)^{a+n_1}} \, I_{(-\infty,r)}(r') \,,
\end{eqnarray}
with $p = \,  n_1 K'/(n_1K' + n_2 J')$, $A_1= \frac{K's + b}{n_2 J'}$, and $A_2= \frac{K's + b}{n_1 K'}$.   Finally, we point out that the above distribution for $R'$ can be described as bilateral Pareto since, conditional on $J',K'$: 
\begin{equation}
\label{paretobilateral}
R'-r|R'>r \,=^{\!\!\!\!^{d}} \hbox{B2}(1, a+n_1-1, A_1)\,, \, r-R'|R'<r \,=^{\!\!\!\!^{d}} \hbox{B2}(1, a+n_1-1, A_2)\, \hbox{ with } p= \mathbb{P}(R'>r)\,.
\end{equation}
Summarizing for the multivariate Lomax case with $K \,=^{\!\!\!\!^{d}} \hbox{G}(\alpha,1)$, the predictive density of $(S',R')$ given $(s,r)$ and prior $\pi$ in (\ref{priorlomax}) may be generated as:
$$ J' \sim h \, , \,  K' \,=^{\!\!\!\!^{d}} \hbox{K}2(\alpha+n_1-1, a+1, \frac{b}{s}, \frac{b}{s}) ,\,$$  and then $R'|J',K'$ and $S'|R',J',K'$ as the bilateral Pareto and Beta type II respectively defined in (\ref{paretobilateral}) and (\ref{Betalomax}) respectively.

\end{example}
\noindent  There are many other potential applications of Theorem \ref{theoremtypeII}, including mixtures of Gamma$(T,\theta)$ distributions with $T \sim g$, but we will focus hereafter on bivariate gamma mixtures, which we define as follows and which is inspired by the Kibble distribution mixture representation (see Example \ref{examplekibble}). 

\begin{definition}
\label{gammabivariatemixture}
A Kibble-type bivariate gamma mixture with parameters $\nu_1, \nu_2, \theta_1, \theta_2$ (all $>0$) and mixing distribution $g$, which we denote $X=(X_1, X_2) \,=^{\!\!\!\!^{d}} \hbox{G}2(\nu_1, \nu_2, \theta_1, \theta_2, g)$, is defined through the representation:
\begin{equation}
\nonumber
X_i|K \, \,=^{\!\!\!\!^{i.d}} \hbox{G}(\nu_i+K, \theta_i), i=1,2, \, \hbox{ with } K \sim g\,.
\end{equation}
\end{definition}

\noindent  In the above definition,  one recovers the original Kibble distribution in (\ref{kibblerepresentation}) for $\nu_1=\nu_2=\nu, \theta_i=\frac{1-\rho}{\lambda_i}$, and
$K \,=^{\!\!\!\!^{d}} \hbox{NB}(\nu, 1-\rho)$, but many other interesting cases are generated including Poisson mixtures with $K \,=^{\!\!\!\!^{d}} \hbox{P}(\lambda)$ and continuous mixing as well.  Of course, one requires $g$ such that $\mathbb{P}(K > - \min{\{\nu_1, \nu_2}\})=1$.  An application of Theorem \ref{theoremtypeII} to such bivariate gamma mixtures yields the following.

\begin{corollary}
\label{corollarygammabivariatemixtures}
Consider type II mixtures as in (\ref{modeltypeII}) with $X|\theta \,=^{\!\!\!\!^{d}} \hbox{G}2(\nu_1, \nu_2, \theta_1, \theta_2, g)$ and $Y|\theta \,=^{\!\!\!\!^{d}} \hbox{G}2(\nu_1', \nu_2', \theta_1, \theta_2, h)$, and prior distribution $\theta_i \,=^{\!\!\!\!^{i.d}} \hbox{G}(a_i, b_i), i=1,2$,  Then,

\begin{enumerate}
\item[{\bf (a)}]   The posterior distribution of $\theta_1, \theta_2$ based on $x$ is 
$\hbox{G}2(\nu_1+a_1, \nu_2+a_2, b_1+x_1, b_2+x_2, g_{\pi,x})$ with $g_{\pi,x}$ a density (or p.m.f.) such that
\begin{equation}
\label{postgammabivariatemxitures}
g_{\pi,x}(k) \, \propto \,  g(k) \,\; \frac{(\nu_1+a_1)_k \, (\nu_2+a_2)_k \, }{(\nu_1)_k (\nu_2)_k} \; \left(  \frac{x_1 \, x_2}{ (b_1+x_1) (b_2+x_2)} \right)^k\,.
\end{equation}

\item[{\bf (b)}]  The predictive distribution for $Y$ admits the representation
\begin{equation}
\label{predgammabivariatemxitures}
Y|J',K' \,=^{\!\!\!\!^{i.d}}  \hbox{B2}(\nu_i'+J', a_i+\nu_i+K', b_i+x_i), i=1,2,   \hbox{ with }
J' \sim h, \, K' \sim g_{\pi,x} \hbox{ independent}.
\end{equation}
\end{enumerate}
\end{corollary}
\noindent {\bf Proof.}  {\bf (a)}  We have $X_i|k' \,=^{\!\!\!\!^{i.d}} \hbox{G}(\nu_i+k', \theta_i)$ and $\theta_i \,=^{\!\!\!\!^{i.d}} \hbox{G}(a_i,b_i)$, which implies that $\pi_{k',x}$ is the density of $\theta_i|k',x \,=^{\!\!\!\!^{i.d}}  \hbox{G}(a_i+\nu_i+k', b_i+x_i)$.   It thus follows from (\ref{posttypeII}) that the posterior distribution of $\theta_1, \theta_2$ is $\hbox{G}2(\nu_1+a_1, \nu_2+a_2, b_1+x_1, b_2+x_2, g_{\pi,x})$.    Finally,  a direct calculation yields
\begin{equation}
m_{\pi, k}(x) \, \propto  \frac{(\nu_1+a_1)_k \, (\nu_2+a_2)_k \, }{(\nu_1)_k (\nu_2)_k} \; \left(  \frac{x_1 \, x_2}{ (b_1+x_1) (b_2+x_2)} \right)^k\,,
\end{equation} 
as a function of $k$, and the result follows.  \\

\noindent {\bf (b)}   This follows from (\ref{preddensitytypeII}), (\ref{postgammabivariatemxitures}), and the evaluation of the predictive density of $Y_i|J=j' \,=^{\!\!\!\!^{i.d}} \hbox{G}(\nu_i'+j', \theta_i)$ based on  $X_i|K=k' \,=^{\!\!\!\!^{i.d}} \hbox{G}(\nu_i+k', \theta_i)$, which has independently distributed $\hbox{B2}$ marginals as given in (\ref{predgammabivariatemxitures}) and previously encountered in Example \ref{exgammaminxtures}.   \qed

\begin{remark}
\label{improperCorollary4.1}
Corollary \ref{corollarygammabivariatemixtures} applies quite generally with 
respect to the type of mixing, and include ``degenerate'' cases with $\mathbb{P}_g(K=0)=1$ and-or  $\mathbb{P}_h(J=0)=1$ corresponding to independent components of $X$ and $Y$ respectively.
The results of Corollary \ref{corollarygammabivariatemixtures} apply as well to some choices of the $a_i$'s and $b_j$'s that lead to an improper prior, as long as the densities $\pi_{K,x}$ exist.  These include for instance the choices $b_1=b_2=0$, with $a_1, a_2$ such that $\mathbb{P}(K> -\max\{\nu_1+a_1, \nu_2+a_2\})=1$.  Observe as well that the choice $a_1=a_2=b_1=b_2=0$, which is typically taken as the usual non-informative prior, leads to the simplification $g_{\pi,x} \equiv g$ in (\ref{postgammabivariatemxitures}) and (\ref{predgammabivariatemxitures}). 
\end{remark}

\begin{example}
\label{exampletypeII}
Corollary \ref{corollarygammabivariatemixtures} applies for any mixing distributions $g$ and $h$, such as Poisson mixing: {\bf (a)} $g \sim \hbox{P}(\lambda_X)$ and $h \sim \hbox{P}(\lambda_Y)$, and negative binomial mixing: {\bf (b)}  $g \sim \hbox{NB}(r_X, \xi_X)$ and $h \sim \hbox{NB}(r_Y, \xi_Y)$, with only the mixing density $g_{\pi,x}$ varying.  We obtain  $$g_{\pi,x} \sim \hbox{Hyp}\left(\nu_1+a_1, \nu_2+a_2; \nu_1, \nu_2; \frac{\lambda_X x_1 x_2}{(b_1+x_1)\,(b_2+x_2)}\right)\,,$$
for the Poisson case {\bf (a)}, and 
$$ g_{\pi,x} \sim \hbox{Hyp}\left(r_X,\nu_1+a_1, \nu_2+a_2; \nu_1, \nu_2; \frac{(1-\xi_X) \, x_1 x_2}{(b_1+x_1)\,(b_2+x_2)}\right)\,,  $$\
for the negative binomial case {\bf (b)}.   Case {\bf (b)} reduces and applies to the original Kibble distribution, as in Example \ref{examplekibble}, for $\nu_1=\nu_2=r_X$, 
$\nu_1'=\nu_2'=r_Y, \xi_X \, = \, \xi_Y \, = \, 1- \rho$, $\theta_i \, = \, \frac{1-\rho}{\lambda_i}, i=1,2$. 
Finally, it is immediate to extend the example to cases where the density $g_{\pi,x}$ is that of a generalized hypergeometric distribution, obtaining
$g_{\pi,x} \sim \hbox{Hyp}(a_1, \ldots, a_p, \nu_1+a_1, \nu_2+a_2; b_1, \ldots, b_q, \nu_1, \nu_2;\frac{\lambda \, x_1 x_2}{(b_1+x_1)\,(b_2+x_2)} )$.
\end{example}

\noindent  Our next example and application of Theorem \ref{theoremtypeII} relates again to Kibble-type bivariate gamma mixtures, but with a prior on $(\theta_1, \theta_2)$ adapted to situations where there exists an order restriction.  Many such situations are plausible.  An interesting case concerns $M|M|1$ queues and the estimation of the arrival ($\lambda$) and service ($\mu$) rates for independent gamma observables, and in particular the ratio $\rho=\lambda/\mu$, whenever a steady-state distribution exists, i.e. $\lambda < \mu$ (Armero \& Bayarri, 1994).  In doing so, it seems natural to consider a McKay distribution as a prior. Another interesting prior choice, which was analyzed by Armero \& Bayarri (1994) for independent gamma observables, is to start with two independent gamma distributions for $\theta_1$ and $\theta_2$ truncated to the ordering set $\theta_2 \geq \theta_1$.   Before presenting how Theorem \ref{theoremtypeII} applies, we catalog some information about a generalized version of the McKay distribution which contains both the above-mentioned prior choices, which will appear below in the posterior analyses, and which is of interest on its own.

\begin{definition} 
\label{definitiongeneralizedmcKay}
The generalized McKay distribution with parameters $  \vec{\gamma}=(\gamma_1, \gamma_2, \gamma_3, \gamma_4, \gamma_5)$, with $\gamma_1>0$, $\gamma_2>0$, $\gamma_4 >0$, $ \gamma_3>\gamma_5 \geq 0$ and $\gamma_1+\gamma_2+\gamma_4>1$, denoted $MG(\vec{\gamma}), $ is a bivariate distribution with p.d.f.  

\begin{equation}
\label{pdfgeneralizedmcKay}
f_{Z_1,Z_2}(z_1, z_2) \, = \, \frac{1}{C( \vec{\gamma})} \, 
\, z_1^{\gamma_1-1} \, (z_2-z_1)^{\gamma_2-1} \, {z_2}^{\gamma_4-1}\, e^{-\gamma_3 z_2 - \gamma_5 z_1} \,  \mathbb{I}_{(0,\infty)}(z_1) \, \mathbb{I}_{(z_1,\infty)}(z_2) \,,
\end{equation}

 \begin{equation}
\label{C}
\hbox{ and } \;  C(\vec{\gamma}) \, = \,  \frac{\Gamma(\gamma_1) \Gamma(\gamma_2)}{\Gamma(\gamma_1 + \gamma_2)} \, \frac{\Gamma(\gamma_1+\gamma_2+\gamma_4-1)}{\gamma_3^{\gamma_1+\gamma_2+\gamma_4-1}}\,_2F_1(\gamma_1, \gamma_1+\gamma_2+\gamma_4-1, \gamma_1+\gamma_2; - \, \frac{\gamma_5}{\gamma_3})\,.
 \end{equation}
\end{definition}
\noindent  Observe that the above distribution reduces to a McKay distribution ($\hbox{M}(\gamma_1, \gamma_2, \gamma_3)$) for $\gamma_4=1$ and $\gamma_5=0$.  As well,  the case $\gamma_2=1$ with $\gamma_5>0$ is generated by two independently distributed $G(\gamma_1, \gamma_5)$ and $G(\gamma_4, \gamma_3)$ truncated to the region $z_2>z_1$.   We have collected some technical observations and properties, including the computation of the normalization constant $C( \vec{\gamma})$, of such distributions in the Appendix.

\begin{example}
Consider the application of Theorem \ref{theoremtypeII} for $X|\theta \,=^{\!\!\!\!^{d}} G2(\nu_1, \nu_2, \theta_1, \theta_2, g)$ and $Y|\theta \,=^{\!\!\!\!^{d}} G2(\nu_1', \nu_2', \theta_1, \theta_2, h)$ combined with the prior $\theta \,=^{\!\!\!\!^{d}} MG(\vec{\gamma})$.   From expressions below for the posterior density of $\theta$ and the predictive density for $Y$, we will be able to extract expressions for Gamma observables with $g$ and $h$ degenerate, as well as for both McKay prior and the truncated Gamma prior considered by Armero \& Bayarri (1994).    Combining the density of $\theta \,=^{\!\!\!\!^{d}} MG(\vec{\gamma})$ with the likelihood$\; f_{\theta,k}(x) \, \propto \, \theta_1^{\nu_1+k} \, e^{-\theta_1 x_1}  \theta_2^{\nu_2+k} \, e^{-\theta_2 x_2}  $, we obtain that
\begin{equation}
\label{pi_kxMG}
\begin{split}
& \qquad \qquad \quad \pi_{k,x}(\theta) \, \sim \, MG(\vec{\gamma}_{k,x}\,), \hbox{ with } \vec{\gamma}_{k,x} =(\gamma_1', \gamma_2', \gamma_3', \gamma_4', \gamma_5'),\\
\, \hbox{ and\ } \gamma_1' & = \gamma_1 + \nu_1 +k, \gamma_2' = \gamma_2, \gamma_3' = x_2 + \gamma_3, \gamma_4' = \gamma_4 + \nu_2 + k, \gamma_5' = x_1 + \gamma_5 \,.
\end{split}
\end{equation} 
We also have 
\begin{eqnarray*}
m_{\pi,k}(x) \, & = &\, \frac{1}{C(\vec{\gamma})} \, \frac{x_1^{\nu_1+k-1} \, x_2^{\nu_2+k-1}}{\Gamma(\nu_1+k)\, \Gamma(\nu_2+k) } \; \int_0^{\infty} \int_{\theta_1} ^{\infty} \theta_1^{\gamma_1'-1} \, (\theta_2 - \theta_1)^{\gamma_2'-1} \, \theta^{\gamma_4'-1} \, e^{-\left(\theta_2 \gamma_3' \, + \, \gamma_5' \theta_1\right)} \, d\theta_2 d\theta_1 \\
\, & = &  \frac{C(\vec{\gamma}_{k,x})}{C(\vec{\gamma})} \, \frac{x_1^{\nu_1+k-1} \, x_2^{\nu_2+k-1}}{\Gamma(\nu_1+k)\, \Gamma(\nu_2+k) }\,, 
\end{eqnarray*}
which implies that 
\begin{equation}
\label{g_pikGM}
g_{\pi,x}(k) \, \propto \,  g(k) \; \frac{(x_1 \, x_2\,)^{k} }{\Gamma(\nu_1+k)\, \Gamma(\nu_2+k) } \,\,  C(\vec{\gamma}_{k,x})\,,
\end{equation}
with $C(\cdot)$ defined in (\ref{C}) and $\vec{\gamma}_{k,x}$ given in (\ref{pi_kxMG}).  From (\ref{posttypeI}), we therefore infer that the posterior distribution $U \,=^{\!\!\!\!^{d}} \pi|x $
admits the representation
\begin{equation}
\label{postMG}
U|K' \,=^{\!\!\!\!^{d}} \hbox{MG}(\vec{\gamma}_{K',x}) \, \hbox{ with }  K' \sim g_{\pi,x}\,,
\end{equation}
as given by (\ref{pi_kxMG}) and (\ref{g_pikGM}).    As mentioned above,  the above representation applies to McKay distributed priors (taking $\gamma_4=1, \gamma_5=0$), as well as the truncated Gamma priors of Armero \& Bayarri (taking $\gamma_2=1, \gamma_5>0$).  However, despite its conciseness, this last representation (\ref{postMG}) remains rather unwieldy due to the $g_{\pi,x}$ density term.  \\

\noindent  Degenerate case of the mixtures $G2$, i.e.,  with $\mathbb{P}(K=0)=1$ in Definition \ref{gammabivariatemixture}, are not particularly interesting when the joint prior distribution for $\theta_1$ and $\theta_2$ factorizes into independent components (e.g., as in Corollary \ref{corollarygammabivariatemixtures}), since the inference problem separates completely into two univariate independent problems.  It is also not admittedly the focus of the mixture models considered in this paper.   However, for non-independent components $\theta_1$ and $\theta_2$ such as the situation here, the Bayesian analysis is also interesting in the degenerate case.  And, the above results applied to 
$X_i \,=^{\!\!\!\!^{i.d}} \hbox{G}(\nu_i, \theta_i)$ yield immediately  the posterior $\theta|x \sim \hbox{MG}(\vec{\gamma}_{0,x})$.  This further connects to the McKay prior by taking $\gamma_4=1$ and $\gamma_5=0$, and to the Armero \& Bayarri prior by taking $\gamma_2=1$ and $\gamma_5>0$.  Furthermore, the latter work considers the change of variables $T_1=\theta_1/\theta_2, T_2=\theta_2$ for reasons of ease of interpretation, and properties of these transformed variables are briefly outlined in the Appendix below. \\

\noindent  Turning now to the Bayesian predictive distribution for $Y|\theta$, still based on  $X|\theta \,=^{\!\!\!\!^{d}} G2(\nu_1, \nu_2, \theta_1, \theta_2, g)$ and prior $\theta \,=^{\!\!\!\!^{d}} MG(\vec{\gamma})$,  we have for (\ref{qpi}), and making use of (\ref{pi_kxMG}),
\begin{eqnarray}
\mkern-18mu \nonumber \ q_{\pi}(y|j',k') \, & = & \,  \int_0^{\infty} \int_{\theta_1}^{\infty}  \prod_{i=1}^2 
\frac{y_i^{\nu_i'+j'-1} e^{-\theta_i y_i} \, }{\Gamma(\nu_i'+j') \, } \,
\, \frac{\theta_i^{(\nu_i'+j')}}{C(\vec{\gamma}_{k',x})} \,
\theta_1^{\gamma_1'-1} \, (\theta_2 - \theta_1)^{\gamma_2'-1} \, \theta^{\gamma_4'-1} \, e^{-\left(\theta_2 \gamma_3' \, + \, \gamma_5' \theta_1\right)}   \;    d\theta_2 d\theta_1 \\
\label{qpiMG}
\, & = & \, \frac{y_1^{\nu_1'+j' -1} \, y_2^{\nu_2'+j'-1} }{\Gamma(\nu_1'+j)\, \Gamma(\nu_2'+j)} \; \frac{C(\vec{\gamma}_{k',j',x,y})}{C(\vec{\gamma}_{k',x})} \, 
\end{eqnarray}
for $y>0$, and with $\vec{\gamma}_{k',j',x,y}\,=\, (\gamma_1'', \gamma_2'', \gamma_3'', \gamma_4'', \gamma_5'') \, , \, \gamma_1''= \gamma_1+\nu_1+\nu_1'+k'+j', \gamma_2''= \gamma_2, \gamma_3''= \gamma_3+x_2+y_2, \gamma_4''= 
\gamma_4+\nu_2+\nu_2'+k'+j'$, and  $\gamma_5''=\gamma_5+x_1+y_1.$  We thus infer that the predictive distribution of $Y$ admits representation (\ref{preddensitytypeII}) with $q_{\pi}(\cdot|J',K')$ given in (\ref{qpiMG}), and $g_{\pi,x}$ given in (\ref{g_pikGM}).  Finally, the same remarks as above applies for the degenerate case with $\mathbb{P}(K=0)= \mathbb{P}(J=0)=1$,  In particular, the Bayesian predictive density for $Y$ is given by 
$q_{\pi}(\cdot|0,0)$, i.e., expression (\ref{qpiMG}) with $j'=k'=0$.

\end{example} 
\noindent  We now turn to a further illustration for models and priors that are bivariate gamma mixtures.   A motivation for this is the flexibility to have prior components that are non-independent in addressing inference for the Kibble distribution parameters or more generally for Kibble type gamma mixtures (Definition \ref{gammabivariatemixture}).  But, in contrast to the previous example, the dependence among the prior components arises not through an order constraint, but by the presence of a hierarchical commonality. 

\begin{corollary}
\label{corollarygammabivariatemixtures+priormixture}
Consider type II mixtures as in (\ref{modeltypeII}) with $X|\theta \,=^{\!\!\!\!^{d}} \hbox{G}2(\nu_1, \nu_2, \theta_1, \theta_2, g)$ and $Y|\theta \,=^{\!\!\!\!^{d}} \hbox{G}2(\nu_1', \nu_2', \theta_1, \theta_2, h)$, and prior distribution $\theta \,=^{\!\!\!\!^{d}} G2(a_1, a_2, b_1, b_2, \epsilon)$.   Then,
\begin{enumerate}
\item[ {\bf (a)}]  The posterior distribution of $\theta$ based on $x$ is $G2(a_1+\nu_1, a_2+\nu_2, b_1+x_1, b_2+x_2, \epsilon_x)$, with the mixing density (or p.m.f.) $\epsilon_x$ is that of $S \, =^d \, K'+L'$ where $(K',L')$ has joint density (or p.m.f.) given by
\begin{equation}
\label{jointK'L'}  p_{K',L'}(k,l) \propto \frac{(a_1+\nu_1)_{k+l} \, (a_2+\nu_2)_{k+l}}{(a_1)_l \, (a_2)_l \, (\nu_1)_k \, (\nu_2)_k}  \, \frac{(x_1x_2)^k \, (b_1 b_2)^l}{\left\lbrace(b_1+x_1)(b_2+x_2) \right\rbrace^{k+l}} \, \;g(k) \,\, \epsilon(l)\,.
\end{equation}

\item[ {\bf (b)}]  The predictive distribution for $Y=(Y_1, Y_2)$ admits the representation
\begin{equation}
\nonumber Y_i|J',K' \,=^{\!\!\!\!^{i.d}}  \hbox{B2}(\nu_i'+J', a_i+\nu_i+K', b_i+x_i), i=1,2,   \hbox{ with }
J' \sim h, \, K' \sim g_{\pi,x} \hbox{ independent},
\end{equation}
\begin{eqnarray} 
\label{g_pi,x(I)}
\hbox{ with } & & \;  g_{\pi,x}(k) \, \propto \,  \frac{ g(k) \, (x_1 x_2)^k}{\{(x_1+b_1)(x_2+b_2)\}^k} \,\,  \frac{I(k,x)}{(\nu_1)_k \, (\nu_2)_k} \,, \\
\nonumber 
 \hbox{ and } & &  I(k,x) \, = \, \int_{\cal{K}}  \frac{(a_1+\nu_1)_{k+l} \, (a_2+\nu_2)_{k+l}}{(a_1)_1 \, (a_2)_l} \;  \left(\frac{b_1 b_2}{(x_1+b_1)(x_2+b_2)} \right)^l  \; d\epsilon(l) \,.
\end{eqnarray}
\end{enumerate}
\end{corollary}
\noindent  {\bf Proof.}
{\bf (a)}  From Theorem \ref{theoremtypeII}, we have 
\begin{equation}
\nonumber
\pi_{k,x}(\theta) \propto \int_0^{\infty} \frac{\theta_1^{a_1+\nu_1+k+l-1} \, \theta_2^{a_2+\nu_2+k+l-1}}{\Gamma(a_1+l) \, \Gamma(a_2+l)} \, e^{-\theta_1(x_1+b_1) - \theta_2(x_2+b_2)} \, (b_1 b_2)^l \, d\epsilon(l)\,,
\end{equation}
which is the density of a $G2\left(a_1+\nu_1+k, a_2+\nu_2+k, x_1+b_1, x_2+b_2, \epsilon'_{x,k}\right)$ distribution having mixing density
\begin{equation}
\nonumber
\epsilon'_{x,k}(l) \, = \frac{1}{I(k,x)} \; \frac{(a_1+\nu_1)_{k+l} \, (a_2+\nu_2)_{k+l}}{(a_1)_1 \, (a_2)_l} \;  \left(\frac{b_1 b_2}{(x_1+b_1)(x_2+b_2)} \right)^l  \; \epsilon(l)\,.
\end{equation}
This along with Theorem \ref{theoremtypeII} tell us that the posterior distribution $U \,=^{\!\!\!\!^{d}} \theta|x$ admits the representation
\begin{equation}
\nonumber   U_i|K',L' \,=^{\!\!\!\!^{i.d}} \hbox{G}(\alpha_i + \nu_i+S, x_i+\beta_i) \,, i=1,2,
\end{equation}
with $S =^d K'+L'$, $L'|K'=k' \sim \epsilon'_{x,k'}$, and $K' \sim g_{\pi,x}$.  Finally, we have by collecting terms and with an interchange in the order of integration: 
\begin{eqnarray}
\nonumber
p_{K',L'}(k,l) & \propto &  \epsilon'_{x,k}(l) \, g(k) \, m_{\pi,k}(x) \\
\nonumber
\,   & \propto &  \epsilon'_{x,k}(l) \, g(k) \, \int_{\Theta} f_{\theta,k}(x) \, \pi(\theta) \, d\tau(\theta) \\
\nonumber
& \propto &  \, \frac{\epsilon'_{x,k}(l) \, g(k) \, (x_1 x_2)^k}{\Gamma(\nu_1+k) \, \Gamma(\nu_2+k)} \, \, \int_{\cal{K}} \int_{\Theta} 
\frac{\theta_1^{a_1+\nu_1+k+l-1} \, \theta_2^{a_2+\nu_2+k+l-1}}{\Gamma(a_1+l) \Gamma(a_2+l)} \, e^{-\theta_1(x_1+b_1) - \theta_2(x_2+b_2)} \, (b_1 b_2)^l \, d\epsilon(l)
\\
\label{g_pi,xII} & \propto & \epsilon'_{x,k}(l) \, g(k) \, \frac{(x_1 x_2)^k}{\{(x_1+b_1)(x_2+b_2)\}^k} \, \frac{I(k,x)}{\Gamma(\nu_1+k) \, \Gamma(\nu_2+k)} \, ,
\end{eqnarray}
which is indeed (\ref{jointK'L'}).

\noindent {\bf (b)}  As in the proof of part (b) of Corollary \ref{corollarygammabivariatemixtures}, the result follows with the evaluation $g_x(k) \propto g(k) \, m_{\pi,k}(x)$ obtained in (\ref{g_pi,xII}) and yielding (\ref{g_pi,x(I)}).   \qed

\begin{remark}
\label{remarkcorollary}
The results of Corollary \ref{corollarygammabivariatemixtures+priormixture} are quite general and many particular cases follow, namely depending on the choices of mixing densities for the model and the prior.    For degenerate mixing under the prior, i.e., $\mathbb{P}_{\epsilon}(L=0)=1$, one recovers Corollary
 \ref{corollarygammabivariatemixtures}.   For independent Gamma distributed observables with $\mathbb{P}_g(K=0)=1$,  Corollary \ref{corollarygammabivariatemixtures+priormixture} applies for the posterior distribution, but with $S=^d L'$ and $L'$ having density $p_{K',L'}(0,l)$ in (\ref{jointK'L'}).  In this last case, Poisson or negative binomial distributed $L'$ leads to results similar to those of Example \ref{exampletypeII}.  Finally, again with $K$ degenerate at $0$, the predictive distribution in part {\bf (b)} of Corollary \ref{corollarygammabivariatemixtures+priormixture} applies with $\mathbb{P}(K'=0)=1$ under $g_{\pi,x}$. 
\end{remark}

\begin{example}
\label{examplekibble-kibble}
We consider Kibble models and priors as a particular application of Corollary \ref{corollarygammabivariatemixtures+priormixture}.   Namely, take $\nu_1=\nu_2$, 
$a_1=a_2$, as well as $K \,=^{\!\!\!\!^{d}} \hbox{NB}(\nu_1, 1-\rho)$, $L \,=^{\!\!\!\!^{d}} \hbox{NB}(a_1, 1-\rho_{\theta})$ as the mixing variables associated with $X|\theta$ and $\theta$ respectively.  With p.m.f.'s  $g(k) \, = \, \frac{(\nu_1)_k}{k!} \, (1-\rho)^{\nu_1} \, \rho^k \,\, \mathbb{I}_{\mathbb{N}}(k)$ and $\epsilon(l) \, = \, \frac{(a_1)_l}{l!} \, (1-\rho)^{a_1} \, \rho^l \,\, \mathbb{I}_{\mathbb{N}}(l)$, we obtain in applying part {\bf (a)} of Corollary \ref{corollarygammabivariatemixtures+priormixture} that
$$ p_{K',L'}(k,l) \, = \, \frac{1}{F_4(a_1+\nu_1, a_1+\nu_1; a_1, \nu_1; \omega_1, \omega_2)} \, \, \frac{(a_1+ \nu_1)_{k+l} \, (a_1+ \nu_1)_{k+l}}{(a_1)_l \, (\nu_1)_k \, l! \, k! }\, \, \omega_1^k \, \omega_2^l\,  \mathbb{I}_{\mathbb{N}}(k) \, \mathbb{I}_{\mathbb{N}}(l)\,,$$ 
with $\omega_1 \, = \, \frac{\rho \, x_1 \, x_2}{(x_1 + b_1)(x_2 + b_2)}$, $\omega_2 \, = \, \frac{\rho_{\theta} \, b_1 \, b_2}{(x_1 + b_1)(x_2 + b_2)}$, 
 and where $F_4$ is the Appell function of the fourth kind given by
$$ F_4(\gamma_1, \gamma_2;  \gamma_3, \gamma_4; w,z) \, = \, 
 \sum_{m=0}^{\infty} \sum_{n=0}^{\infty} \, \frac{(\gamma_1)_{m+n} \, (\gamma_2)_{m+n} \, }{m! \, n! \, (\gamma_3)_n \, (\gamma_4)_m \, } \, w^m \, z^n\,.$$
Part {\bf (b)} of Corollary \ref{corollarygammabivariatemixtures+priormixture} applies directly for the Bayesian predictive distribution with $J' \,=^{\!\!\!\!^{d}} \hbox{NB}(a_1, 1-\rho_{\theta})$ and $g_{\pi,x}$ derived from (\ref{g_pi,x(I)}) or equivalently as the marginal p.m.f. of $K'$ for the $p_{K',L'}$, which is given by
\begin{equation}
g_{\pi,x}(k) \, = \, \frac{_2F_1(a_1+\nu_1+k, a_1+\nu_1+k;a_1; \omega_2)} {F_4(a_1+\nu_1, a_1+\nu_1; a_1, \nu_1; \omega_1, \omega_2)} \, \, \, \frac{(a_1+ \nu_1)_{k} \, (a_1+ \nu_1)_{k}}{\, (\nu_1)_k \,  k! }\, \, \omega_1^k \,  \mathbb{I}_{\mathbb{N}}(k) \,.
\end{equation}
\end{example}

\begin{remark}
\label{expectation-posterior}
The posterior distribution in Corollary \ref{corollarygammabivariatemixtures} is conveniently expressed in terms of a $\hbox{G}2$ mixture distribution with mixing variable $S=K'+L'$ based on joint density (or p.m.f.) given in (\ref{jointK'L'}).   The previous example gives a more precise form for a negative binomial mixture.   The properties of $S$ remain to be explored in general or for other specific choices of $g$ and $\epsilon$.   In terms of the posterior expectation, it follows readily from the above that:
\begin{eqnarray*}
\mathbb{E}(\theta|x) \,  =  \, \mathbb{E}^{S|x} \left\lbrace \mathbb{E}(\theta|x,S)  \right\rbrace  \, \,  =   \left(\frac{a_1+\nu_1+ \mathbb{E}(S|x)}{x_1+b_1}  \,, \,
\frac{a_2+\nu_2+ \mathbb{E}(S|x)}{x_2+b_2} \,
           \right) \,.
\end{eqnarray*} 
In the context of Example \ref{examplekibble-kibble},  one obtains with standard manipulations that 
$$  \mathbb{E}(S|x) \, = \, \mathbb{E}(K') \, + \, \mathbb{E}(L') \,= \, \frac{(a_1 \, + \, \nu_1)^2 \, (\omega_1 A_2/\nu_1 + \omega_2 A_3/a_1)}{\, A_1\,} \, \,, $$
with  $A_1\, = \, F_4(a_1+\nu_1, a_1+\nu_1; a_1, \nu_1; \omega_1, \omega_2)$, $A_2\, = \, F_4(a_1+\nu_1 +1, a_1+\nu_1+1; a_1, \nu_1+1; \omega_1, \omega_2)$, and $A_3\, = \, F_4(a_1+\nu_1+1, a_1+\nu_1+1; a_1+1, \nu_1; \omega_1, \omega_2)$. 
\end{remark}

\section*{Appendix}
\noindent  Here are some observations on the generalized McKay distribution $MG(\vec{\gamma}), $ with p.d.f. given in (\ref{pdfgeneralizedmcKay}).

\begin{enumerate}
\item[ {\bf A.}]    With the transformation $(Z_1, Z_2) \to (T_1=Z_1/Z_2, T_2=Z_2)$, one obtains the density
\begin{equation}
\label{densityT1T2}
f_{T_1, T_2}(t_1, t_2) \, = \, \frac{1}{C(\vec{\gamma})} \,   t_1^{\gamma_1 - 1} \, (1-t_1)^{\gamma_2-1}\, t_2^{\gamma_1+\gamma_2+\gamma_4-2} \, e^{-t_2(\gamma_3+\gamma_5 t_1)} \,\,   \mathbb{I}_{(0,1)}(t_1) \, \mathbb{I}_{(0,\infty)}(t_2)\,.
\end{equation}
From this, it follows immediately that 
$$ T_2|T_1 \,=^{\!\!\!\!^{d}} \hbox{G}(\gamma_1+\gamma_2+\gamma_4-1, \gamma_3+\gamma_5 T_1) \,,$$
and then using the value of $C(\vec{\gamma})$ given in (\ref{C}), we obtain the marginal density 
\begin{equation}
\nonumber
f_{T_1}(t_1) \, = \,  \frac{\Gamma(\gamma_1+\gamma_2)}{\Gamma(\gamma_1) \, \Gamma(\gamma_2)} \, \frac{t_1^{\gamma_1 - 1} \, (1-t_1)^{\gamma_2-1}}{_2F_1(\gamma_1, \gamma_1+\gamma_2+\gamma_4-1, \gamma_1+\gamma_2; -\frac{\gamma_5}{\gamma_3})}
\, \left(\frac{\gamma_3}{\gamma_3+\gamma_5 t_1}\right)^{\gamma_1+\gamma_2+\gamma_4-1} \,
 \, \mathbb{I}_{(0,1)}(t_1)\,.
\end{equation}
Observe that we have independence between $T_1$ and $T_2$ iff $\gamma_5=0$, which includes the McKay case, and that $T_1 \,=^{\!\!\!\!^{d}} \hbox{B}(\gamma_1, \gamma_2)$ whenever $\gamma_5=0$. We point out that, for $\gamma_4=1$, the above density for $T_1$ matches a generalized Beta distribution given by Chen \& Novick (1984) with parameters $\gamma_1, \gamma_2, \lambda = 1+ \frac{\gamma_5}{\gamma_3}$, and p.d.f.  $  \lambda^{\gamma_1} \, \frac{\Gamma(\gamma_1+\gamma_2)}{\Gamma(\gamma_1) \Gamma(\gamma_2)} \, \frac{t_1^{\gamma_1 - 1} \, (1-t_1)^{\gamma_2-1}}{\left(1 -(1-\lambda) t_1 \right)^{\gamma_1+\gamma_2}} \, \mathbb{I}_{(0,1)}(t_1) $.  More generally, the density has appeared in Exton (1976), as well as in Armero \& Bayarri (1994) where it is referred to as a Gauss hypergeometric distribution.

\item[ \bf B.]   The given expression for the normalization constant $C(\vec{\gamma})$ can be checked by integrating (\ref{densityT1T2}) via an expansion of the term $e^{-\gamma_5 t_1 t_2}$ in its MacLaurin series, and interchanging integrals and sum.  

\item[ {\bf C.}]   The marginals densities of $Z$, which are gamma distributed for the McKay distribution, are given as follows:
\begin{eqnarray*}
f_{Z_1}(z_1) \, & =  \, & \frac{1}{C(\vec{\gamma})}  \int_{z_1}^{\infty}   f_{Z_1,Z_2}(z_1, z_2) \, dz_2 \\    \, & =  \, & \frac{1}{C(\vec{\gamma})} \, z_1^{\gamma_1 - 1} \, e^{-\gamma_5 z_1} \, \int_{0}^{\infty} \, u^{\gamma_2 -1} \, (u+z_1)^{\gamma_4 - 1} \, e^{-\gamma_3 (z_1 + u)} \, du \\ 
\, & =  \, & \frac{\Gamma(\gamma_2)}{C(\vec{\gamma})} \, z_1^{\gamma_1 \, + \, \gamma_2 \, + \gamma_4 - 1} \, e^{-(\gamma_5+\gamma_3) \, z_1} \;\psi(\gamma_2, \gamma_2+\gamma_4, \gamma_3 \, z_1)\,,
\end{eqnarray*}  
for $z_1>0$, using the identity $\int_{0}^{\infty} t^{a-1} \, (\Delta+t)^{b-a-1} \, e^{-ct} \, dt \, = \, \Delta^{b-1} \, \Gamma(a) \, \psi(a,b,c\Delta)$ which holds for $a,c>0$ and $b \in \mathbb{R}$. Observe that $Z_1 \,=^{\!\!\!\!^{d}} \hbox{G}(\gamma_1, \gamma_3+\gamma_5)$ for $\gamma_4=1$, as $\psi(a,a+1, c) \, = \, c^{-a}$ for $a,c>0$, and that the additional specification $\gamma_5=0$ leads us back to the McKay case. \\

\noindent  For the marginal density of $Z_2$, we have
\begin{eqnarray*}
f_{Z_2}(z_2) \, & =  \, & \frac{1}{C(\vec{\gamma})} \, \,   z_2^{\gamma_4-1} \, e^{-\gamma_3 z_2} \, \int_{0}^{z_2} \, z_1^{\gamma_1-1} \, (z_2-z_1)^{\gamma_2-1} \, e^{-\gamma_5 z_1} \, dz_1 \\  \, & = \, &  \frac{1}{C(\vec{\gamma})} \, \,   z_2^{\gamma_1+\gamma_2+ \gamma_4-1} \, e^{-\gamma_3 z_2} \int_{0}^{1}  \, u^{\gamma_1-1} \, (1-u)^{\gamma_2-1} \, e^{-\gamma_5 z_2 u } \, du \\  \, & = \, & \frac{z_2^{\gamma_1+\gamma_2+ \gamma_4-1} \, e^{-\gamma_3 z_2}}{\Gamma(\gamma_1+\gamma_2+\gamma_4-1)} \,  \, \frac{\gamma_3^{\gamma_1+\gamma_2+\gamma_4-1} \, _1F_1\left(\gamma_1; \gamma_1+\gamma_2; -\gamma_5 z_2\right)}{_2F_1\left(\gamma_1, \gamma_1+\gamma_2+\gamma_4-1; \gamma_1+\gamma_2; -\frac{\gamma_5}{\gamma_3} \right)}\,,
\end{eqnarray*}
for $z_2>0$, and with the change of variables $z_1= z_2 u$, as well as the confluent hypergeometric identity  $\int_0^1 u^{a-1} (1-u)^{b-1} \, e^{-c u} \, du \, = \,  \frac{\Gamma(a) \Gamma(b)}{\Gamma(a+b)} \, _1F_1(a; a+b;-c)$, $a,b,c>0$.  Observe that $Z_2 \,=^{\!\!\!\!^{d}} \hbox{G}(\gamma_1+\gamma_2+\gamma_4-1, \gamma_3)$ for $\gamma_5=0$, while
 $f_{Z_2}(z_2) \, = \, \frac{(\gamma_3^{\gamma_2}) \, (\gamma_3+\gamma_5)^{\gamma_1}}{\Gamma(\gamma_1+\gamma_2)} \, 
z_2^{\gamma_1+\gamma_2} \, e^{-\gamma_3 z_2} \, _1F_1(\gamma_1, \gamma_1+\gamma_2; -
\gamma_5 z_2) \, \mathbb{I}_{(0,\infty)}(z_2)$ for $\gamma_4=1$ using the identity 
$_1F_0(a; -; -c) \, = \, (1+c)^{-a}$ for $a>0$ and $c \in (-1,1)$.

\item[ {\bf D.}]  The moments are readily expressible in terms of the constant of normalization $C(\vec{\gamma})$.  For instance, it is easy to obtain the mixed moments   
$$  \mathbb{E}(Z_1^m \, Z_2^n\,) \, = \, \frac{(\gamma_1)_m \, (\gamma_2)_n \, (\gamma_1+\gamma_2+\gamma_4-1)_{s}}{(\gamma_1+\gamma_2)_{s} \, \gamma_3^{s}} \;  
\frac{_2F_1(\gamma_1+m, \gamma_1+\gamma_2+\gamma_4+s-1, \gamma_1+\gamma_2+s; -\frac{\gamma_5}{\gamma_3})}{_2F_1(\gamma_1, \gamma_1+\gamma_2+\gamma_4-1, \gamma_1+\gamma_2; -\frac{\gamma_5}{\gamma_3})}\,,$$
for $m,n \in \mathbb{N}$ and $s=m+n$.

\item[ {\bf E.}] For McKay distributed $(Z_1, Z_2)$, the variables $U_1=Z_1$ and $U_2=Z_2-Z_1$ are independently distributed.   For the general case, it is easy to verify that we obtain conditional distributions that are Kummer type II, namely:
$$ U_2|U_1 \,=^{\!\!\!\!^{d}}  \hbox{K}2(\gamma_2 \, , \,  1-\gamma_2-\gamma_4 \, , \, \gamma_3 \, U_1 \,, \, U_1) \hbox{ and }  U_1|U_2 \,=^{\!\!\!\!^{d}}  \hbox{K}2(\gamma_1 \, , \,  1-\gamma_2-\gamma_4 \, , \, (\gamma_3+\gamma_5) \, U_2 \,, \, U_2)\,,$$   
with independence iff $\gamma_4=0$.
\end{enumerate}

\section*{Concluding remarks} 
\noindent For two common types of mixtures, we have provided analytical, novel and practical representations of Bayesian posterior and predictive distributions corresponding to a given prior distribution.   By doing so, we have added to the catalogue of such Bayesian solutions which require little computation except recourse to some special but familiar functions and distributions.   Specific model and prior distribution choices have been put to the forefront for illustrative purposes and for their practicality.  These include inference for common non-central distributions (chi-square, Beta, Fisher), the distribution of the coefficient of determination $R^2$ in multiple regression models, mixtures of normal and exponential distributions, and bivariate gamma mixtures including the Kibble distribution.   

\section*{Acknowledgements}

\noindent \'Eric Marchand gratefully acknowledges the research support from the 
Natural Sciences and Engineering Research Council of Canada.

\end{document}